\documentclass[11pt]{article}
\usepackage{amsthm,amssymb,amsmath}
\usepackage{bm}
\usepackage{soul}
\usepackage{subeqnarray}
\usepackage[top=1.5cm, bottom=3cm, left=2.5cm, right=2.5cm]{geometry}
\usepackage[utf8]{inputenc}
\usepackage{graphicx}
\usepackage{booktabs}
\usepackage{cases}
\usepackage{caption}
\usepackage{setspace}
\usepackage{color}
\usepackage{xcolor}
\usepackage{subfigure}
\usepackage{comment}
\newcommand{\cS}{{\cal S}} 
\newcommand{\cV}{{\cal V}} 
\newcommand{\cT}{{\cal T}} 
\newcommand{\bU}{{\bf U}} 
\newcommand{\cX}{{\cal X}} 
\newcommand{\cY}{{\cal Y}} 
\newcommand{\cP}{{\cal P}} 
\newcommand{\cQ}{{\cal Q}}

\numberwithin{equation}{section}

\begin{document}


\title{Dynamic tipping and cyclic folds, in a one-dimensional non-smooth dynamical system linked to climate models.}

\author{Chris Budd \thanks{Dept. of Mathematical Sciences, University of Bath, UK} \and Rachel Kuske \thanks{School of Mathematics, Georgia Institute of Technology, Atlanta, USA.}} 

\maketitle

\begin{abstract}{We study the behaviour at tipping points close to (smoothed) non-smooth fold  bifurcations in one-dimensional oscillatory forced systems. 
The focus is the Stommel-Box, and related climate models, which are piecewise-smooth continuous dynamical systems, modelling thermohaline circulation. These exhibit non-smooth fold bifurcations which arise when a saddle-point and a focus meet at a border collision bifurcation. By using techniques from the theory of non-smooth dynamical systems we are able to provide precise estimates for the general tipping behaviour at the non-smooth fold as parameters vary. These are significantly different from the usual tipping point estimates, showing a much more rapid rate of tipping. We also see very rapid, and non-monotone, changes in the tipping points due to the effect of non-smoothness in the system. All of this has important implications for the prediction of tipping in climate systems.}
 \end{abstract}


\section{Introduction}

\subsection{Overview}

\noindent Tipping behaviour plays an important role in many dynamical systems, particularly those arising as models of climate. See \cite{Lenton11} for a number of different cases of dynamical systems where tipping arises. It is hard to define tipping precisely; however, essentially it is the phenomenon of large scale changes in the dynamics of a system (such as a major change in the climate) which arise when a critical parameter $\mu$ is varied in some manner. Such a variation could be a slow drift in $\mu$, a periodic excitation (such as a seasonal climatic variation), or even a stochastic perturbation. In 
\cite{Lenton11,paillard2004antarctic,stommel1961thermohaline} a series of models of climatic states where tipping behaviour is thought to be likely are given, ranging from the loss of the Arctic sea ice, to changes in the mid-Atlantic over turning circulation. Associated with these models are estimates of when tipping will occur, and the consequent impact on our climate. Further studies of the detection of tipping points (using either statistical or machine learning methods) are described in \cite{Lenton11}, \cite{bury}. A common feature of all of the climate models considered above is that they assume that the dynamical systems under study are {\em smooth}. In such cases there is a fairly complete understanding of when, and how, tipping can occur \cite{zhu2015tipping}. 
However, this smoothness assumption is only partly valid in the context of models of climate. The vastly differing timescales in climate models, and the possibility of thresholding or similar behaviour means that, to leading order, many climate models are non-smooth, and have differential equations with right hand sides that are non-smooth functions of their arguments. Three examples are the PP04 model of glacial cycle \cite{MoruBudd20} which takes the form of a Filippov system) \cite{bernardo2008piecewise}, the Budyko-Sellars model of ice-albedo feedback, and the Stommel 2-box model for thermohaline circulation in North Atlantic which we study in this paper \cite{kaper},\cite{stommel1961thermohaline}. The latter system takes the form of a system of differential equations  with  nonlinearities which lose differentiability for certain values of their arguments. Non-smooth systems have additional, {\em discontinuity induced} bifurcations, to those studied in smooth systems. An important such example being a {\em non-smooth fold} (NSF) where two fixed points coalesce in a 'V-shaped' bifurcation diagram, in contrast to the well-known parabolic shape near a classical smooth fold or saddle node bifurcation (SNB).  Our interest in this paper lies in the climatically important question of tipping behaviour close to a NSF. 

\vspace{0.1in}

\noindent In contrast to the case of smooth systems, the study of tipping in non-smooth dynamical systems is much less developed. The case of static bifurcation at a NSF is described analytically in \cite{bernardo2008piecewise}. Furthermore in \cite{KimWang18} a series of numerical results are given for stochastically induced tipping close to a NSF arising in an electronic system. Calculations of tipping close to the NSF in the Stommel 2-box model are presented in \cite{kaper}, \cite{dijkstra2013nonlinear}. However, a complete analytical study of these systems, particularly in the context of climate dynamics, and hence the impact on the understanding of the dynamics, and detection, of tipping close to a NSF, remains open.

\subsection{The non-smooth model studied}

\noindent In this paper we study tipping in a one-dimensional reduction of the Stommel 2-box model, which whilst simple enough to analyse carefully, has enough complexity to exhibit significantly different behaviour from the classical tipping at a saddle-node bifurcation (SNB) studied in the dynamical systems literature. In particular we study the piece-wise linear system

\begin{equation}
    \frac{dx}{dt} =  2|x| - \mu(t) + f(t), \quad x(0) = x_0.
    \label{c10}
\end{equation}
Here $x$ is the state variable, $f(t)$ is a (seasonal) forcing, and $\mu(t)$ a (climatic) control variable, for which changes lead to tipping in $x$. For a derivation of this non-smooth one-dimensional system from the Stommel 2-box model, see \cite{kaper}.  The system (\ref{c10}) was studied in the earlier paper \cite{cody} using multi-scale methods of analysis. In this paper we significantly extend the results in \cite{cody} by using methods from the theory of non-smooth dynamical systems. In doing this we make extensive use of the fact that the piece-wise linear nature of the nonlinearity allows us to construct the solution as a sequence of exact solutions to a linear ODE. Whilst the piece-wise linear function $|x|$ is of course an approximation to the nonlinear models $F(x,\mu)$ encountered in climate dynamics (for example the Budyko-Sellars or Stommel 2-box models), we are justified in considering (\ref{c10}) with the nonlinearity $2|x|$ to be a useful normal form for a NSF for a general nonlinear function $F(x,\mu)$ if the non-smooth fold point of $F$ occurs when $x = \mu = 0$ and $|x|$ is small \cite{bernardo2008piecewise}. Of course when tipping occurs and $|x|$ becomes large, then this approximation will break down. Our assumption then is that the system evolves into some different asymptotic state, for example into a stable fixed point of an extended nonlinear system with $dx/dt = F(x,\mu)$ (for example in the Stommel 2-box model studied in \cite{cody}). This latter dynamics is not of interest in this paper which concentrates on tipping behaviour only.

\vspace{0.1in}

\noindent Trivially, if $f = 0$ and $\mu > 0$ is constant, then (\ref{c10}) has a {\em stable} fixed point at $x^- = -\mu/2$ and an {\em unstable} fixed point at $x^+ = \mu/2$. There is a static (non-smooth fold NSF) at 
\begin{equation}
    \mu \equiv \mu_{NSF} = 0,
    \label{feb1}
    \end{equation}
where the fixed points $x^-$ and $x^+$ coalesce at $x = 0$. 
If $\mu < 0$ then the system has no fixed points, and $x \to \infty$ in infinite time, with 
$$x(t) \sim C e^{2t} \quad \mbox{as} \quad t \to \infty.$$
Hence we see a profound change in the behaviour of this system for positive and negative $\mu$. For dynamically varying $\mu$ and non-zero $f$ we see rather more complex behaviour. The change in the control parameter $\mu(t)$  is given in this paper by taking
\begin{equation}
\mu(t) = \mu_0 - \epsilon t, \quad 0 < \epsilon \ll 1 
\label{c12}
\end{equation}
representing slow climatic change. Similarly the function $£f(t)$  is given by 
\begin{equation}
    f(t) = A \cos(\omega t).
    \label{c11.5}
\end{equation}
representing a more rapid periodic seasonal perturbation. 
\vspace{0.1in}

\noindent Our scenario for tipping is that as $\mu$ decreases from a positive value, then tipping occurs for some critical value of $\mu_{TP}.$ We aim to estimate this value as a function of $A, \epsilon$ and $\omega$, using a combination of non-smooth dynamical systems methods and asymptotic analysis. 

\subsection{A 'non-smooth' definition of tipping}\label{tipdef}

\noindent The solutions of the non-smooth system (\ref{c10}) and \eqref{c12}  behave in an {\em essentially different way} from those where tipping occurs close to a saddle-node bifurcation (SNB). The latter problem has been studied extensively in the literature, (see for example \cite{Lenton11}), and is given by
\begin{eqnarray} \frac{dx}{dt} = \frac{x^2}{\alpha} - \mu(t). \label{SNBS}
\end{eqnarray}
In this case the solution diverges to infinity in a finite time, and tipping arises when $\mu(t) = \mu_{TP}$ where
\begin{eqnarray}
\mu_{TP}  \sim - c_0 \;  \alpha^{1/3} \epsilon^{2/3},  \label{1.7}
\end{eqnarray}
which we note is lagged relative to the static SNB at $\mu=0$.
Here $c_0 = 2.3381 \ldots$ is the first zero of the reversed Airy function ${\rm Ai}(-x)$.

\vspace{0.1in}

\noindent In contrast the non-smooth system (\ref{c10}) that we study, is {\em linear} for large $x > 0$ and the solution diverges to infinity exponentially in an {\em infinite time}. In this case we define the non-smooth problem (\ref{c10}) to have {\em tipped} at a time $t_{TP}$, with $\mu_{TP} = \mu_0 - \epsilon t_{TP}$ when 
\begin{equation}
x(t_{TP}) = K > 0, \quad \mbox{and} \quad \dot{x}(t_{TP}) > 0,
\label{Kcond}
\end{equation}
for some sufficiently {\em large} value $K$. Using this definition we show in general that the system 
with a NSF tips earlier than the SNB case, directly related to smaller asymptotic scale ${\cal O}(\epsilon \log\epsilon)$ of lagged tipping for a varying control parameter through the static NSF bifurcation, in contrast to ${\cal O}(\epsilon^{2/3} ) $   for the SNB \eqref{1.7} . Furthermore,  in both settings  external forcing $f(t)$ contributes to an advance of the tipping $\mu_{TP}$,  which increases with A and variation of  $\omega$. This advance due to external forcing may dominate over the lag due to varying $\mu(t)$, leading to tipping earlier than the static bifurcation,  $\mu_{TP} > \mu_{NSF} = 0$. The precise relationship between the tipping point $\mu_{TP}$, the slow drift rate of $\mu$, and the strength and frequency of the seasonal forcing is subtle,  with the possibility of large variations in the value of $\mu_{TP}$ as these parameters vary. It is this relationship that we explore in this paper.

\subsection{Results}

\vspace{0.1in}

\noindent We will look at the following three cases of (\ref{c10}) 

\begin{enumerate}
    \item {\em Slow drift.} In this case we take $f(t) = 0$ and set $d\mu/dt = -\epsilon < 0.$ The solution in this case has a quasi-steady state which tips at a critical value of $\mu_{TP} < \mu_{NSF}$.
    
    \item {\em Oscillatory only forcing.} In this case we take $f(t) = A \cos(\omega t)$ and $\mu$ fixed. For small $A$ or sufficiently large $\mu$ we see a stable periodic solution. This exists for all $\mu > \mu_{CF}$ and ceases to exist at a cyclic fold $\mu_{CF}$. The system is unstable if $\mu < \mu_{CF}$ and we see tipping in this case. We study this system for both large and small values of $\omega$, and give precise estimates for $\mu_{CF}$ in both cases. 
    \item {\em Oscillatory forcing and slow drift.} In this case we let $\epsilon > 0$ and take $f(t) = A \cos(\omega t)$. This tips at a value of $\mu_{TP} < \mu_{CF}.$ For large $\omega$ the value of $\mu_{TP}$ is monotonic in $\omega$ or $\epsilon$ and can be computed using the multi-scale arguments in \cite{cody}. For smaller values of $\omega$ the value of $\mu_{TP}$ is not monotonic, and we see sharp gradients in its dependence on $\omega$.
    
\end{enumerate}

\vspace{0.1in}

\noindent For example, amongst other results we will establish the following.

\vspace{0.2in}

\noindent {\bf Lemma A} (Following \cite{cody}) {\em If $A = 0$, $\omega = 0$,  and 
$d\mu/dt = -\epsilon$ with $\quad 0 < \epsilon \ll 1,$ 
then tipping in the non-smooth system (\ref{c10}) arises when, to leading order, 
$$\mu_{TP} \sim -\epsilon \log(K/\epsilon)/2.$$
}
We contrast this with the estimate for tipping at the SNB given in (\ref{1.7}).

\vspace{0.2in}

\noindent {\bf Lemma B} {\em (i) If $A > 0, \omega \ne 0, \epsilon = 0$ then problem (\ref{c11}) has a stable periodic solution. If $\mu < \mu_G = A/\sqrt{1 + \omega^2/4}$ this takes both positive and negative values. This solution exists down to a cyclic fold when $\mu = \mu_{CF} < \mu_G$. 
For $\mu < \mu_{CF}$ the periodic solution ceases to exist and we see a divergence of the solution to infinity. 

\vspace{0.1in}

\noindent (ii) For large $\omega$ we have
$$\mu_{CF} \sim \frac{4A}{\pi \omega} \left( 1 - \frac{L}{\omega^2} \right),$$
where $L \approx 0.7 \ldots$,
and for small $\omega$ we have
$$\mu_{CF} \sim \frac{A}{1 + \omega^2/4}.$$

}

\vspace{0.2in}

\noindent {\bf Lemma C} {\em If $A > 0, \omega \ne 0, \epsilon > 0 $ then problem (\ref{c11}) has a tipping point for some $\mu_{TP} < \mu_{CF} < \mu_G.$ The value of $\mu_{TP}$ decreases monotonically for large $\omega$, but as either $\omega$ or $\epsilon$ decrease to zero the we see large transitions in its value.

}

\vspace{0.1in}

\begin{figure}[!htbp]
	\centering
		  \includegraphics[width=3.in]{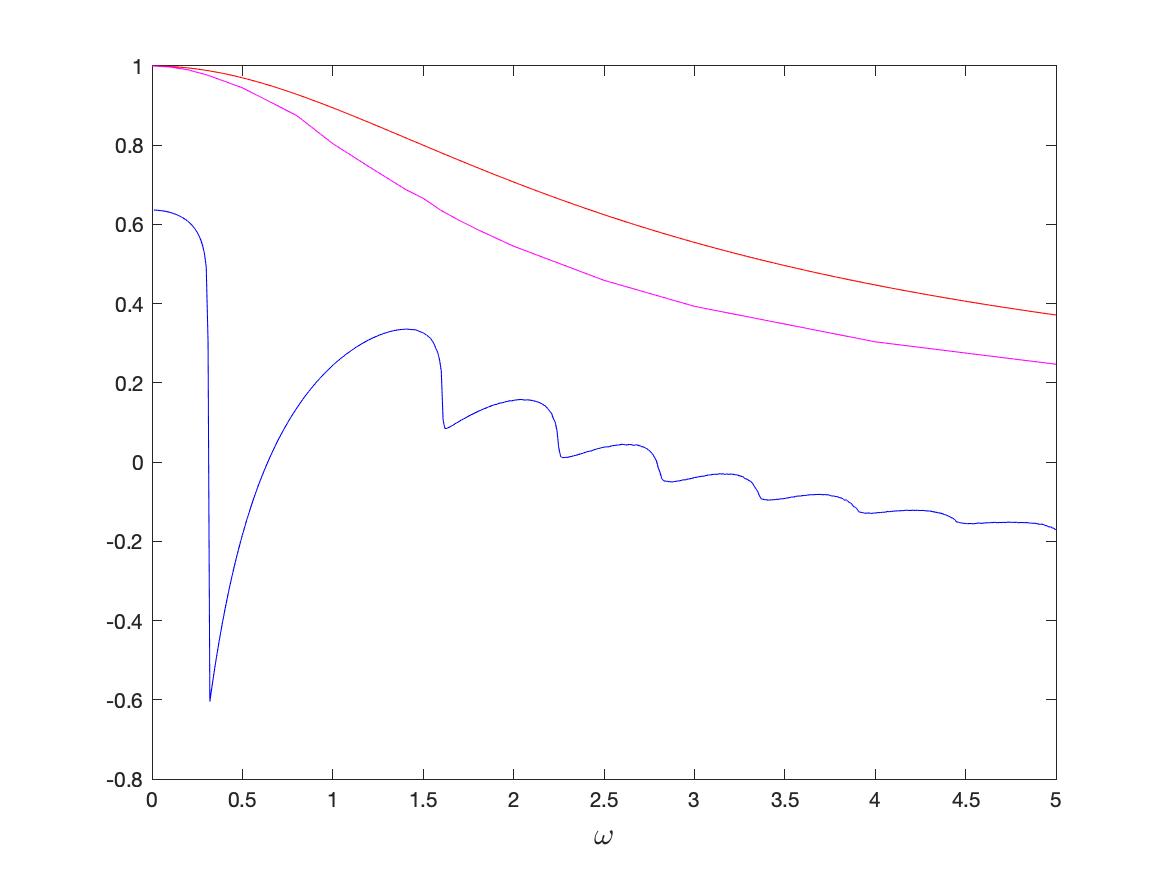}
    \includegraphics[width=3.in]{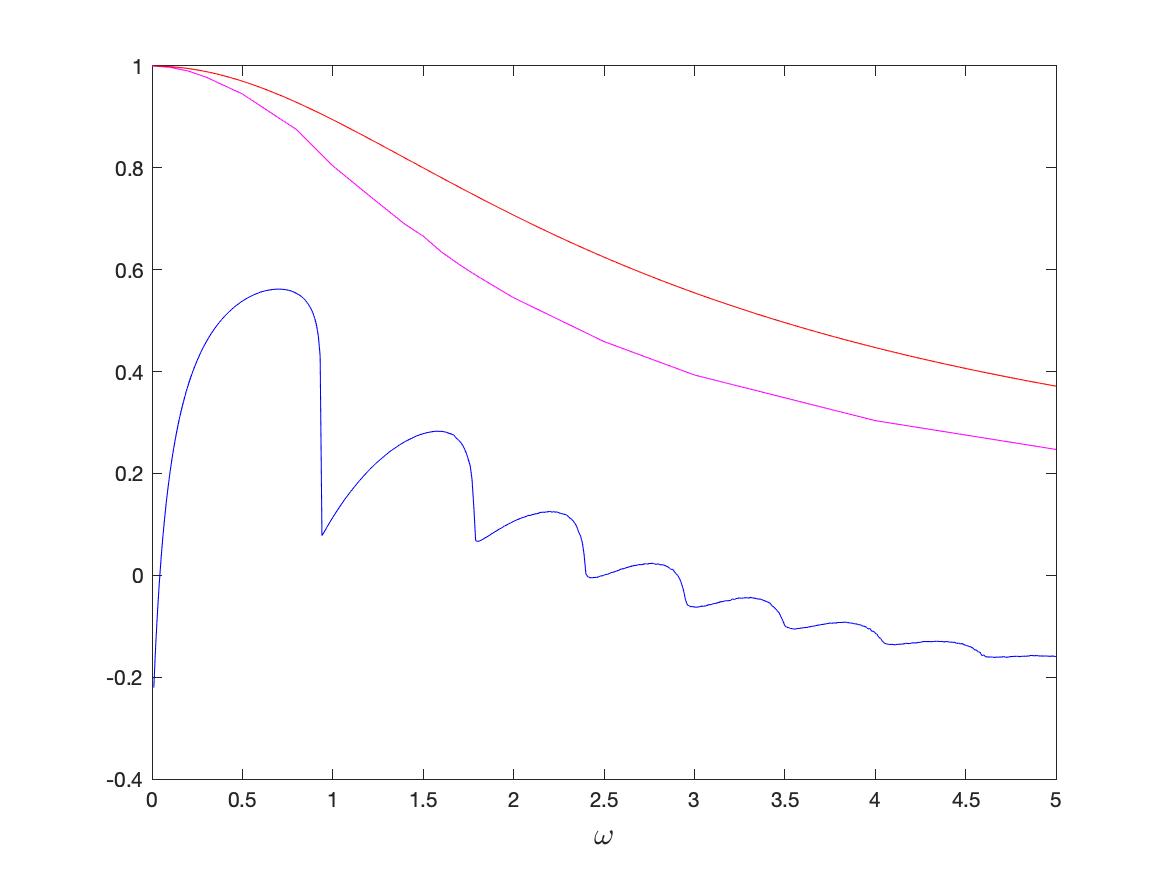}
	\caption{Calculated values of $\mu_{TP}$ (blue), $\mu_{CF}$ (maroon) and $\mu_G$ (red) as a function of $\omega$, when $\epsilon = 0.1, K = 10$ and $A = 1$. Here we take $\mu = 1 - \epsilon t$ and $x(0) = -1/2$. 
 On the left we take forcing $A \cos(\omega t)$ and on the right forcing $A \sin(\omega t)$. Note that in both cases we may have $\mu_{TP} > \mu_{NSF}.$} 
	\label{fig:dec1}
\end{figure}

\noindent We illustrate the conclusions of Lemma C in Figure \ref{fig:dec1} (left) where we fix the drift rate $\epsilon = 0.1$, set $A = 1, K = 10$, and vary $\omega$. For initial conditions we take $\mu(0) = 1$ and $x(0) = -1/2$ (which is the point $x^-$ when $\mu = 1$). In this figure we have plotted the tipping point $\mu_{TP}$ in blue, the cyclic fold $\mu_{CF}$ in maroon and the grazing point $\mu_G$ in red. Note that the curve for $\mu_{TP}$ whilst regular for large $\omega$ as predicted in \cite{cody}, has a complex form for smaller $\omega$, with a sharp transition at $\omega \approx 0.32$ and further transitions for larger $\omega$. As a comparison in Figure \ref{fig:dec1} (right) we plot the same figure, but with the forcing $A \sin (\omega t)$ in (\ref{c10}). This latter figure shows the impact of the phase of the forcing, on the tipping. This will be discussed in Section 4. 

\vspace{0.2in}

\noindent {\em Smoothing} Non-smooth systems typically arise as limits of smooth systems as certain parameters tend to zero;  for example, in a problem with two different time-scales, one of which is much shorter than the other. It is therefore meaningful to compare problem (\ref{c10}) in the broader context of a system of smoother models. A smoother version of (1.1) is given by
\begin{equation}
    \frac{dx}{dt} = 2 \sqrt{\alpha^2 + x^2} - 2\alpha - \mu(t) + f(t), \quad x(0) = x_0,
    \label{c11}
\end{equation}
for some $0 \ll \alpha \ll 1.$
The interplay between the smoothing $\alpha$ and the other parameters, such as the tipping rate $\epsilon$, is subtle. If $\alpha \ll \epsilon $ then the dynamics resembles that of the non-smooth system (\ref{c10}), however if $\alpha \gg \epsilon$ then the behaviour is closer to that seen near to a saddle-node bifurcation (\ref{1.7}). In all cases we find that $\mu_{TP}$ {\em decreases} as the smoothing $\alpha$ {\em increases}. In other words, {\em smoothing postpones tipping.}

\vspace{0.1in}

\noindent {\em Implications} The climatic implications of the above results are significant. Systems which have non-smooth (or close to non-smooth) governing equations are likely to having tipping earlier than that predicted by the classical analysis. But the point at which tipping occurs may depend in a non-monotonic manner on the various tipping parameters. These results and the sharp transitions in $\mu_{TP}$ have important implications in climate studies and lead to a level of uncertainty in determining when tipping will occur in a general system. 

\vspace{0.1in}

\subsection{Summary}

\noindent The remainder of this paper is structured as follows. In Section 2 we briefly study tipping in  the non-smooth system (\ref{c10}) under slow drift without forcing. In Section 3 we  use methods from the theory of non-smooth dynamical system to study the {\em cyclic folds} $\mu_{CF}$ in the periodic solutions of (\ref{c10}) with no drift and periodic forcing. We look at the limits of both small and large $\omega$. In Section 4 we combine the results of Sections 2 and 3 to look at tipping, complementing the previous results of \cite{cody}. In the case of large $\omega$ we find that $\mu_{TP}$ decays smoothly with $\omega$. However, as $\omega$ decreases then we find that $\mu_{TP}$ experiences large transitions. We also consider the impact of the initial conditions, and also the phase of the forcing, on the location of the tipping point. In Section 5 we partially extend all of these results to the smoothed system (\ref{c11}) where we see that smoothing postpones tipping. We give some  further analytic results for $A=0$ combined with some numerical calculations for $A\neq 0$. 
Finally in Section 6 we draw some climate related conclusions from this work. Proofs of various technical results are given in the  appendices.

\section{Tipping with slow drift only}

\noindent We first consider the system (\ref{c10}) without forcing, and with slow drift given by
$d\mu/dt = -\epsilon.$ This problem was studied in \cite{cody} and from this we have the following result:

\vspace{0.2in}

\noindent {\bf Lemma 1} {\em Take $x(0) = x_0 < 0$ and $d\mu/dt = -\epsilon$, $\mu(0) \equiv \mu_0 = {\cal O}(1)$  (in general, $\mu_0$ is an $O((1)$ distance from $\mu_{NSF}$). Then if $\epsilon \ll 1$, the tipping value  for slow drift, defined as $\mu_{\epsilon}$ is given asymptotically for $\epsilon \ll 1$ by
\begin{equation}
    \mu_\epsilon = -\frac{\epsilon}{2} - \frac{\epsilon}{2} \log \left(\frac{2K}{\epsilon}\right) - \frac{\epsilon^2}{8K} \log \left( \frac{2K}{\epsilon} \right)  + \ldots.
    \label{ca7}
\end{equation}
}

\vspace{0.1in}

\noindent {\em Proof } The proof makes use of the piece-wise linearity of the function $|x|$. This allows us to solve the equation (\ref{c10}) exactly in the two regions $x < 0$ and $x > 0$. We can then match these two solutions to get the result. Technical details of this construction are given in Appendix A, and in \cite{cody} for similar models.

\section{Cyclic Folds in the oscillatory forced system without drift}

\subsection{Overview}

\noindent We next consider (\ref{c10}) in the case of zero drift $\epsilon = 0$ and with oscillatory forcing. Our emphasis will be on the region of existence of the periodic solutions of the resulting system. In \cite{cody} these were studied using a multi-scale averaging method, assuming large $\omega$. Here we take a different approach and instead look at the solution of the algebraic equations satisfied by the periodic solution.
This approach is applicable for the full range of values of $\omega$, for which we are able to investigate the bifurcation behaviour of these solutions by using methods from the theory of non-smooth dynamical systems. In particular we look at the periodically forced system:
\begin{equation}
    \frac{dx}{dt} = 2|x| - \mu + A \cos(\omega t). 
    \label{o1}
\end{equation}

\noindent In this system we can regard any of the parameters $\omega$, $A$, and $\mu$  as appropriate bifurcation parameters. We note that if we set
$x = A y, \quad \mu = A \nu$
then we have
\begin{equation}
    \frac{dy}{dt} = 2|y| - \nu + \cos(\omega t).
    \label{o2}
\end{equation}
so we may always, without loss of generality, set $A = 1$ and then rescale. We will use an algebraic method to look at both the large, and the small, $\omega$ limits of the solutions.

\vspace{0.1in}

\noindent {\em General solution behaviour} If $x < 0$ for all $t$ then there is a stable periodic solution to (\ref{o1}) of period $T = 2\pi/\omega$ of the form
\begin{equation}
x(t) = -\frac{\mu}{2} + \frac{ A}{\sqrt{4 + \omega^2}} \cos(\omega t - \phi). 
\label{n1}
\end{equation}
This periodic solution always exists in the region $x < 0$ if $\mu$ is sufficiently large. If $\omega$ and $A$ are {\em fixed}, and $\mu$ decreases, then there is a {\em grazing event} at $\mu = \mu_G$ at which  the periodic solution in (\ref{n1})  grazes the discontinuity manifold $\Sigma$ so that there is a time $t$ at which 
$x(t)=x'(t)=0, x''(t) <0$. This occurs when 

\begin{equation}
\mu \equiv \mu_G = \frac{A}{\sqrt{1 + \omega^2/4}}.
\label{n3}
\end{equation}
Note that $\mu_G \sim 2A/\omega$ as $\omega \to \infty.$

\vspace{0.1in}

\noindent As $\mu$ is reduced further then there is still a periodic solution of (\ref{o1}). There are sub-intervals $t \in [0,2\pi/\omega] $, when $x(t) > 0,$ and other sub-intervals when $x < 0$.
The periodic solution persists for $\mu < \mu_{G}$ until it loses stability at a smooth {\em Cyclic Fold bifurcation} at $\mu = \mu_{CF} < \mu_G.$ 




\vspace{0.1in}

\noindent If there is no drift, then the system diverges to infinity if $\mu < \mu_{CF}$ and is stable if $\mu > \mu_{CF}$. 

\vspace{0.1in}

\noindent In this section we obtain asymptotic estimates for $\mu_{CF}$ for both large and small values of $\omega$.  We find from our numerical calculations of $\mu_{CF}$ that these two estimates have a region where they overlap. Hence we may approximate $\mu_{CF}$ over the full range of $\omega$. These estimate improve on the results of  \cite{cody} where a multi-scale (averaging) analysis of (\ref{o1}) for large $\omega$ gave the leading order estimate

\begin{equation}
\mu_{CF} \sim \frac{4A}{\pi \omega}.
\label{n10}
\end{equation}

\noindent We state our new estimates as follows.

\vspace{0.2in}

\noindent {\bf Lemma 2} {\em (i) If $\omega$ is {\em large} there is a constant $L$ such that

\begin{equation}
\mu_{CF} \sim \frac{4 A}{\pi \omega } \left(1 - \frac{L}{\omega^2} \right) \quad \mbox{as} \quad \omega \to \infty.
\label{2}
\end{equation}

\vspace{0.1in}

\noindent (ii) In contrast, for small values of $\omega$ we have

\begin{equation}
\mu_{CF} \sim \frac{A}{1 + \omega^2/4} \quad \mbox{as} 
\quad \omega \to 0.
\label{2.2}
\end{equation}

}

\noindent {\bf NOTES}

\begin{enumerate}
    \item  The asymptotic theory we have developed does not give the value of $L$ directly. However extensive numerical experiments indicate that a value of $L \approx 0.7 \ldots$ fits the data to high accuracy.

\item The two estimates (\ref{2}) and (\ref{2.2}) overlap when $\omega \approx 2$. We find that each gives a very good description of the location of the cyclic fold in the respective ranges of $\omega > 2$ and $\omega < 2$. In all cases the numerically computed curve of the values of $\mu_{CF}$ lies close to the maximum of the two estimates.

\end{enumerate}

\subsection{Set up of the algebraic equations}

\noindent We assume that $\mu_{CF} < \mu < \mu_G$, and $A$, are fixed and that we  look at a non-smooth solution of (\ref{c10}) which intersects the manifold $\Sigma: x = 0,$ and spends a non-zero time in the region $S^+ : x > 0$  and a non-zero time in the region $S^-: x < 0.$ We  also assume that this solution is periodic and synchronised to the forcing, so that it has period 
$$T = 2 \pi n/\omega.$$

\noindent NOTE There is so far no numerical evidence of {\em sub-harmonic solutions} for which $n > 1$. We will take $n = 1$ throughout. 

\vspace{0.1in}

\noindent To construct the solution we  assume that the solution $x(t) \equiv x^+(t)$ lies in $S^+$ for a time interval $a/\omega < t < b/\omega$ and the solution $x(t) \equiv x^-(t)$ lies in $S^-$ for a time interval $b/\omega < t < a/\omega + 2\pi /\omega.$ Determining the   values of $a$ and $b$ is part of the asymptotic analysis.
We then have the (local) {\em compatibility conditions}

\begin{equation}
    x^+(a/\omega) = x^+(b/\omega) = x^-(b/\omega) = x^-(a/\omega + 2\pi/\omega) = 0.
    \label{4}
\end{equation}

\noindent In each interval we can exploit the fact that the system is piece-wise linear to write down an {\em exact solution}. In particular, 
there are constants $C^+$ and $C^-$ so that:
\begin{equation}
x^+(t) = \frac{\mu}{2} + C^+ e^{2(t - a/\omega)} + Q^+(t), \qquad 
x^-(t) = -\frac{\mu}{2} + C^- e^{2(t - b/\omega)} + Q^-(t),
\label{v1}
\end{equation}
where the periodic functions $Q^{\pm}(t)$ are given by:
\begin{equation}
Q^{\pm}(t) = \frac{A}{4 + \omega^2} \left( \omega \sin(\omega t) \mp 2\cos(\omega t) \right). 
\label{v2}
\end{equation}


\noindent Substituting this expression into the compatibility conditions (\ref{4}) leads to  four algebraic equations satisfied by four unknowns $a,b,C^+,C^-$.

\begin{equation}
    \frac{\mu}{2} + C^+  + Q^+(a) = 0,
    \label{7}
    \end{equation}
   
\begin{equation}
    \frac{\mu}{2} + C^+e^{2(b-a)/\omega}  + Q^+(b) = 0,
    \label{8}
    \end{equation}

\begin{equation}
    -\frac{\mu}{2} + C^-  + Q^-(b) = 0,
    \label{9}
    \end{equation}
    
\begin{equation}
    -\frac{\mu}{2} + C^- e^{-2(2\pi + a - b)/\omega} + Q^-(a) = 0.
    \label{10}
    \end{equation}    
     

\vspace{0.1in}

\begin{figure}[htb!]
    \centering
    \includegraphics[width=0.45\textwidth]{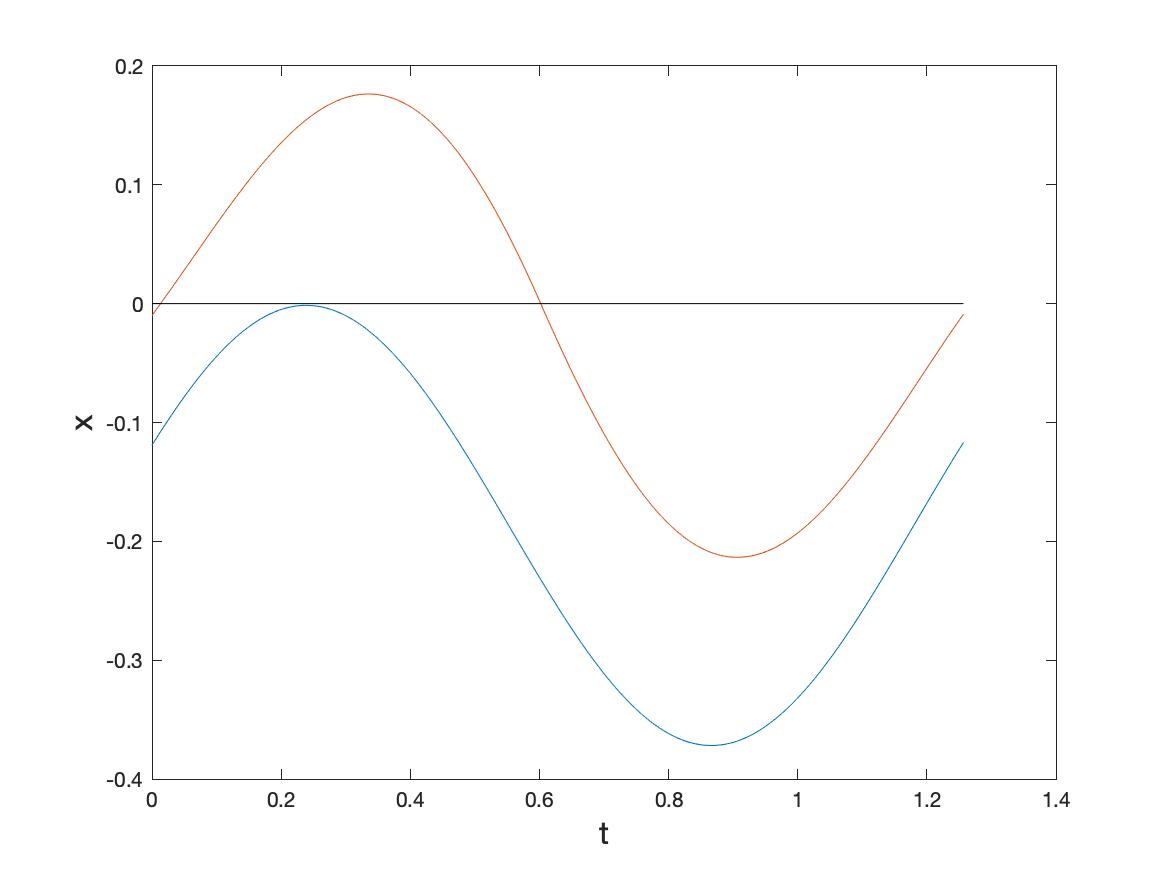}
    \includegraphics[width=0.45\textwidth]{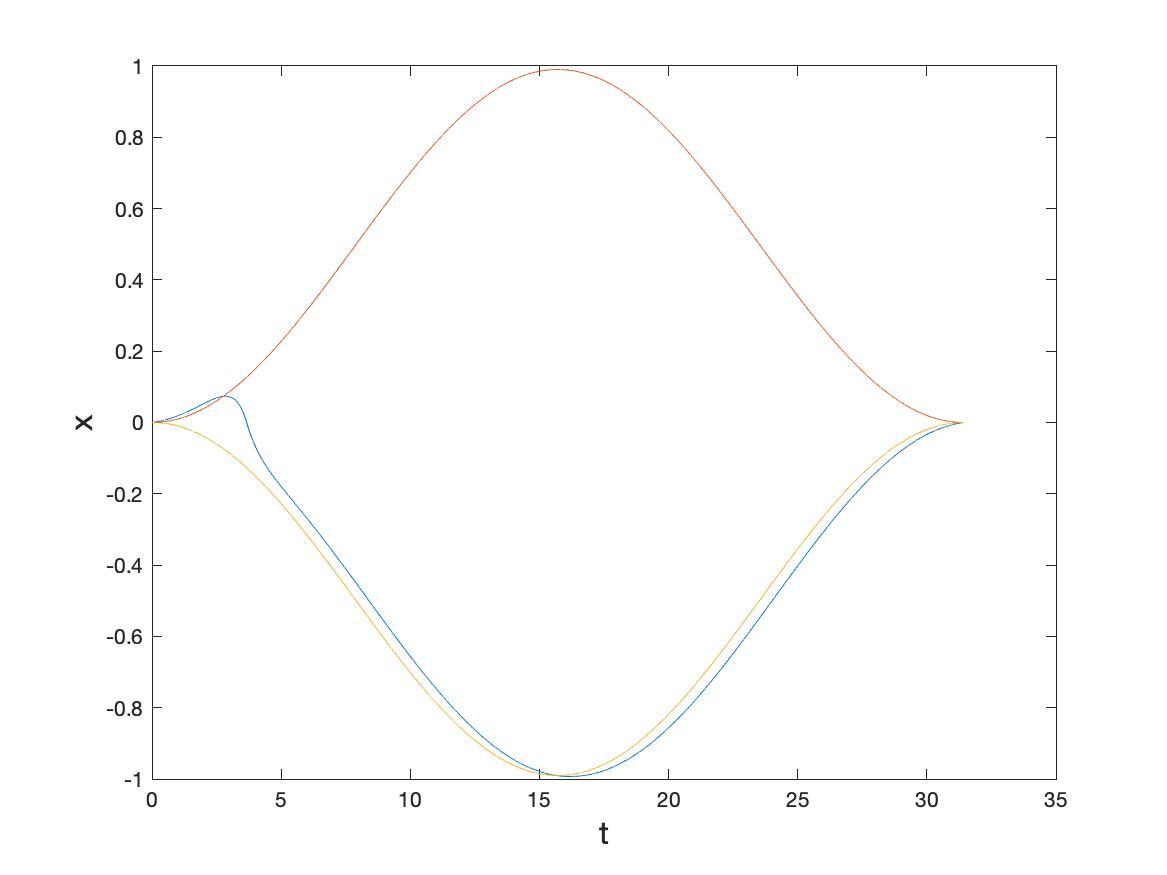}
    \caption{(left) The solution $x(t)$ when $\omega = 5$ is large. Below is the solution at the grazing value $\mu_G =0.3714$ with the graze clearly visible. Above is the solution close to the cyclic fold $\mu_{CF} = 0.2471$. Note that for $\mu$ near $\mu_{CF}$ the solution is close to being symmetric about zero, with a small zero point $a$ and small mean $\langle x \rangle$. 
    (right) The solution $x(t)$ (in blue) close to the cyclic fold when $\omega = 0.2$ is small, and $\mu = 0.9902$, $A = 1$. Also shown are the curves $\pm \mu (1 - \cos(\omega t))/2.$ Note that the solution only spends a small proportion of the time in the region $x > 0$.}
    \label{fig:2novaaa}
\end{figure}

\noindent The system (\ref{7})-(\ref{10}) gives a set of equations for the four unknowns $C^+, C^-, a,b.$ We can consider $\mu$ to be a bifurcation parameter, with the solutions depending continuously on $\mu$. These conditions apply only to the solutions which exist when $\mu_{CF} < \mu < \mu_{G}$. If $\mu > \mu_{G}$  we have $x^-$ only, as given in \eqref{n1} so the compatibility conditions are not relevant, and we have 
$C^- = 0$ at $\mu = \mu_G$. 
We now proceed to estimate $\mu_{CF}$ for both large and small values of $\omega$. To help motivate this calculation we show in Figure 2 the forms that $x(t)$ can take in various limits.
If $\omega = 5$ and $A = 1$ then $\mu_G = 0.3714$ and $\mu_{CF} =0.2471$. A plot of $x(t)$ for these values of $\mu$ is given in the left pane of Figure \ref{fig:2novaaa}, with the solution for $\mu = \mu_G$ taking values below those for $\mu = \mu_{CF}$. The graze is clearly visible on the lower figure and the solution has mean $\langle x \rangle = -\mu_G/2$. When $\mu = \mu_{CF}$ we see that $x(t)$  is close to being symmetric about the line $x = 0$, with both the zero point $a$, and the mean $\langle x \rangle$, close to zero. We shall establish these results rigorously in Lemma 3. 
In contrast, in the right pane of  Figure \ref{fig:2novaaa} we plot the solution at the cyclic fold when $\omega = 0.2$ is small. The two panels of Figure 2 illustrate the fundamental difference between large and small $\omega$ that manifests itself in the estimates for $\mu_{CF}$, namely symmetric behavior near $\mu_{CF}$ for large $\omega$, and asymmetric behavior for small $\omega$.

\subsection{Estimates of $\mu_{CF}$ for large $\omega$ with $A$ fixed}

\noindent In \cite{cody} a multi-scale averaging approach for $\omega \gg 1$ is used to give the estimate (\ref{n10}) for $\mu_{CF}$. We now use an asymptotic approach applied to the algebraic system (\ref{7})-(\ref{10}) to find a more precise large-$\omega$ approximation to $\mu_{CF}$, with the main asymptotic result given in (\ref{2}). We do this by obtaining estimates for $a$,  $C^+, C^-$ and the mean $C$ of $x(t)$ when $\omega \gg 1$ and $A$ is fixed, looking at the cases of general $\mu$ and $\mu = \mu_{CF}$.  These are given by:

\vspace{0.2in}

\noindent {\bf Lemma 3} {\em Let $A = {\cal O}(1)$. Then as $\omega \to \infty$

\vspace{0.1in}

\noindent (i) If $\mu_{CF} < \mu < \mu_{G}$  then $a = {\cal O}(1)$.

\vspace{0.1in}

\noindent (ii) If $\mu = \mu_{CF}$ then $a = {\cal O}(1/\omega )$.

\vspace{0.1in}

\noindent (iii) For all $\mu_{CF} \le \mu \le \mu_G$
\begin{eqnarray}
C^- - C^+ = \mu + {\cal O}(a/\omega^2), \quad C^-+C^+ = {\cal O}(a/\omega), \quad 
\langle x \rangle = {\cal O}(a/\omega).
\end{eqnarray}
}

\noindent These results help to explain the behaviour shown in the upper curve in the left pane of Figure \ref{fig:2novaaa} for which  $\mu = \mu_{CF}$. In this case $a \sim {\cal O}(1/\omega)$, and hence the mean $\langle x \rangle \sim {\cal O}(1/\omega^2) $, so that both become small when $\mu$ is close to the cyclic fold.  In contrast, the lower curve in this pane shows the solution when $\mu = \mu_G$. In this case neither $a$ nor the mean $\langle x \rangle$ are close to zero.

\vspace{0.1in} 

\noindent Finally, as a key to obtaining these results, and predicting the form of the solution close to the cyclic fold bifurcation when $\omega$ is large, we have the following estimate for $a$ in terms of $\mu$. 

\vspace{0.2in}

\noindent {\bf Lemma 4} If $\omega \gg 1$ and $A = {\cal O}(1)$ is fixed, then
\begin{equation}
    \pi \mu = \frac{4A}{\omega} \left[ a \sin(a) + \cos(a) \right] + {\cal O}\left( \frac{a}{\omega^2}, \frac{1}{\omega^3} \right).
    \label{coct1}
\end{equation}

\begin{figure}[htb!]
    \centering
    \includegraphics[width=0.45\textwidth]{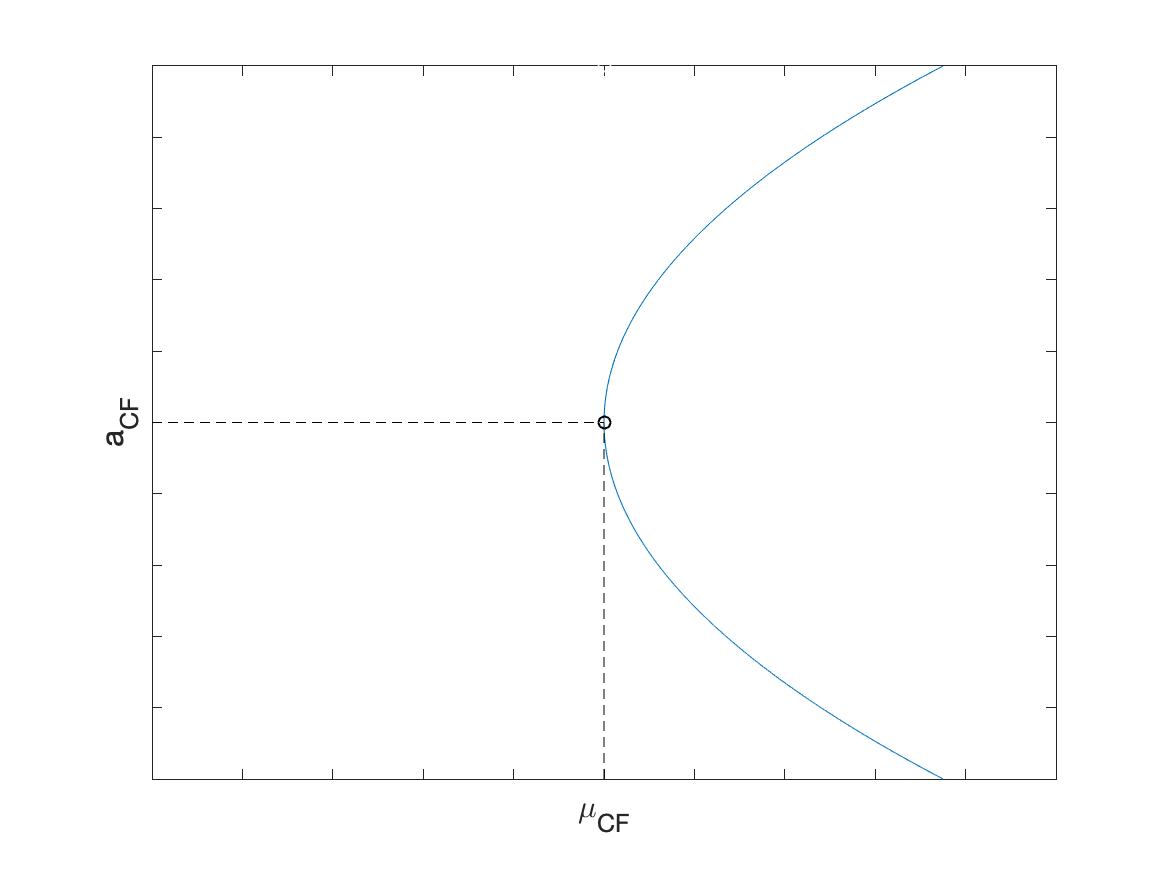}
    \caption{The solution $(\mu,a)$ of (\ref{coct1}) showing the cyclic fold at $(\mu_{CF},a_{CF})$ with $\mu_{CF} \sim 4A/(\pi \omega).$}
    \label{cf1}
\end{figure}

\noindent A schematic plot of the solution $a$ of (\ref{coct1}) in terms of $\mu$ is given in Figure \ref{cf1}. The estimate of $a$ and $\mu_{CF}$ given by solving (\ref{coct1}) then leads to the estimates in Lemma 3, and then to the result in Lemma 2 (i). 

\noindent We see that this cyclic fold exhibits the classical shape of the Saddle-Node bifurcation (SNB), with a quadratic form close to $\mu_{CF}$.

\vspace{0.1in}

\noindent {\em Proof} The derivation of Lemma 2(i), and of Lemmas 3 and 4 is technical, and is given in Appendix B.

\subsection{Estimates of $\mu_{CF}$ for small $\omega$}

\noindent We observe from the above that for small $\omega$ we have $\mu_{CF} \approx \mu_G$ (with $\mu_{CF} < \mu_G$). This is because, as $\omega \ll 1$ the solution can only spend a small time in the region $x > 0$ due to the strong exponential growth in this region. 
As a simple first calculation we set $s = \omega t$ giving
$$ \omega \frac{dx}{ds} = 2|x| - \mu + A \cos(s).$$
Hence to leading order we have
$$2|x| =  \mu - A \cos(s).$$
We will show that as $\omega \to 0$ then $A/\mu = 1 + {\cal O}(\omega^2)$ and hence
\begin{equation}
2|x| = \mu (1 - \cos(s)).
\label{nov2}
\end{equation}

\noindent The full solution $x(t)$ is shown in Figure \ref{fig:2novaaa} (right) for the case of $\mu = 0.9902, A = 1$ and $\omega = 0.2$. In this figure we can see the solution crossing from $\mu(1 - \cos(\omega t))/2$ to $-\mu(1 - \cos(\omega t))/2$. The crossing occurs relatively close (expressed as a fraction of its period) to zero as $\omega \to 0.$ This picture gives us a lot of insight into the general solution behaviour. However, it is very difficult to continue the  perturbation analysis directly in the differential equation as in \eqref{nov2}, due to the non-smoothness of the system. Accordingly we now look directly at the algebraic system (\ref{7}-\ref{10})  in the case of $\omega \ll 1$.

\vspace{0.1in}

\noindent Noting that as $\omega \to 0$ that $\exp(-2\pi/\omega)$ is zero to all orders, we may (to all orders) write 
(\ref{7}) and (\ref{10}) as:

\begin{equation}
\frac{\mu}{2} + C^+ + \frac{A}{4 + \omega^2} \left[ \omega \sin(a) - 2 \cos(a) \right] \equiv \frac{\mu}{2} + C^+ - \frac{A}{\sqrt{4 + \omega^2}} \cos(a + \phi) = 0,
\label{nov3}
\end{equation}

\begin{equation}
-\frac{\mu}{2} + \frac{A}{4 + \omega^2} \left[ \omega \sin(a) + 2 \cos(a) \right] \equiv -\frac{\mu}{2} + \frac{A}{\sqrt{4 + \omega^2}} \cos(a - \phi)= 0.
\label{nov4}
\end{equation}
From (\ref{nov4}) we have
\begin{equation}
a = \phi + \arccos \left( \frac{\mu \sqrt{1 + \omega^2/4}}{A} \right), \quad \tan(\phi) = \omega/2.
\label{nov5}
\end{equation}
Adding and subtracting (\ref{nov3},\ref{nov4}) gives, respectively,
\begin{eqnarray}
C^+ = -\frac{2 A \omega \sin(a)}{4 + \omega^2} = {\cal O}(\omega)  \qquad  C^+ = \frac{4A \cos(a)}{4 + \omega^2} - \mu.
\label{nov5}
\end{eqnarray}

\noindent We observe that the coefficient $C^+$ corresponds to {\em exponentially increasing} terms in the region for $x > 0$. If $C^+ > 0$ these will grow rapidly, and the solution will not cross the line $x = 0.$ Hence we must have $C^+ < 0$ and hence $a > 0$. A sufficient condition for this to occur is that
\begin{equation}
A < (1 + \omega^2/4) \; \mu.
\label{nov7}
\end{equation}

\noindent Now we consider the equation (\ref{8}). The most important term in this expression for our analysis is given by $C^+ e^{2(b-a)/\omega}.$
We note that if $b-a = {\cal O}(1)$  then this expression is exponentially unbounded as $\omega \to 0$. Hence, in order to obtain a bounded solution in this limit we require that there is an ${\cal O}(1)$ quantity $\Delta$ such that 
\begin{equation}
\Delta = \frac{b-a}{\omega}, \quad 
b - a = \omega \Delta.
\label{nov88}
\end{equation}
We then write equation (\ref{8}) as
$$\frac{\mu}{2} + C^+ e^{2\Delta} = \frac{A}{\sqrt{4 + \omega^2}} \cos(a + \phi + \omega \Delta), $${
Substituting in (\ref{nov4},\ref{nov5}), rearranging,  and dividing through by $A$ we have
$$\cos(a - \phi) = \cos(a + \phi + \omega \Delta) + \frac{2 \omega \sin(a) e^{2 \Delta}}{\sqrt{\omega^2 + 4}}.$$

\vspace{0.1in}

\noindent We note first that as $\omega \to 0$ we have
$a = {\cal O}(\omega).$ To see this, suppose not. Then expanding the above to ${\cal O}(\omega)$ and noting that  from (\ref{nov5}) $\phi \sim \omega/2$, we have at ${\cal O}(\omega)$: 
$$\sin(a)/2 = -\sin(a) (1/2 + \Delta) + \sin(a) e^{2 \Delta}.$$
If $\sin(a) \ne 0$ then dividing by $\sin(a)$ and adding, gives:
$1 + \Delta = e^{2 \Delta}$
which is a contradiction. Hence to this order we have $\sin(a)=0$ and hence $a=0$. Accordingly we set $a = \omega a_1$ to give
$$\cos(\omega(a_1 - 1/2)) = \cos(\omega(a_1 + 1/2 + \Delta)) + \frac{2 \omega \sin(\omega a_1) e^{2 \Delta}}{\sqrt{\omega^2 + 4}}.$$
Expanding to ${\cal O}(\omega^2)$, and rearranging, we find that
\begin{equation}
a_1 = \frac{\Delta^2 + \Delta}{2 \left( e^{2 \Delta} - 1 - \Delta \right) }.
\label{nov9}
\end{equation}

\vspace{0.1in}

\noindent We deduce that in the limit as $\omega \to 0$ that in the solution of (\ref{7}-\ref{10}) we have

$$\mu_{CF} = \frac{A}{1 + \omega^2/4}, \quad a  = \omega a_1 > 0, \quad b - a = \omega \Delta. $$

\vspace{0.1in}

\noindent From (\ref{nov9}) we can estimate $a_1$ in terms of $\Delta$, noting that as $\Delta$ increases then $a_1 \sim \Delta^2 \exp(-2\Delta)/2$, and hence decreases very rapidly.

\vspace{0.1in}

\noindent This concludes the proof of Lemma 2 (ii).

\qed 


\subsection{Numerical estimates of the cyclic fold}
\label{numerics}

\noindent {\em Overview}: To give numerical support for the above calculations we fix $A = 1$ and solve the algebraic system (\ref{7}-\ref{10}) directly using the Matlab routine {\tt fsolve}, combined with a path-following method, for a variety of values of $\omega$ and $\mu$.  To do this it is easiest to fix $\omega$ and to start the calculation when $\mu = \mu_G$ at which $C^- = 0$. We then take $\mu$ to be a bifurcation parameter, and slowly reduce it until the solver {\tt fsolve} fails at $\mu = \mu_{CF}$.
This is done for the whole range of values of $\omega \in [0,20]$. The results are given in the table in Appendix B.

\vspace{0.1in}

\noindent {\em Large $\omega$}: The values of $\mu_{CF}$ for $\omega > 2$ are shown in Figure \ref{fig:dec20}. The asymptotic theory (\ref{2}) predicts that to order $1/\omega^2$ we have
$$\mu_{CF} \sim \frac{4A}{\pi \omega} \left(1 - \frac{L}{\omega^2} \right)$$
where the constant $L$ is to be determined. Accordingly we plot
$( 1 - \pi \omega \mu_{CF}/4A ) $ as a function of $\omega$ on a log-log graph, and compare it with a plot of $0.7/\omega^2$ on the same graph. We see from this that the asymptotic form is correct, and we make an estimate of $L = 0.7$ 
\begin{figure}[htb!]
    \centering
    \includegraphics[width=0.4\textwidth]{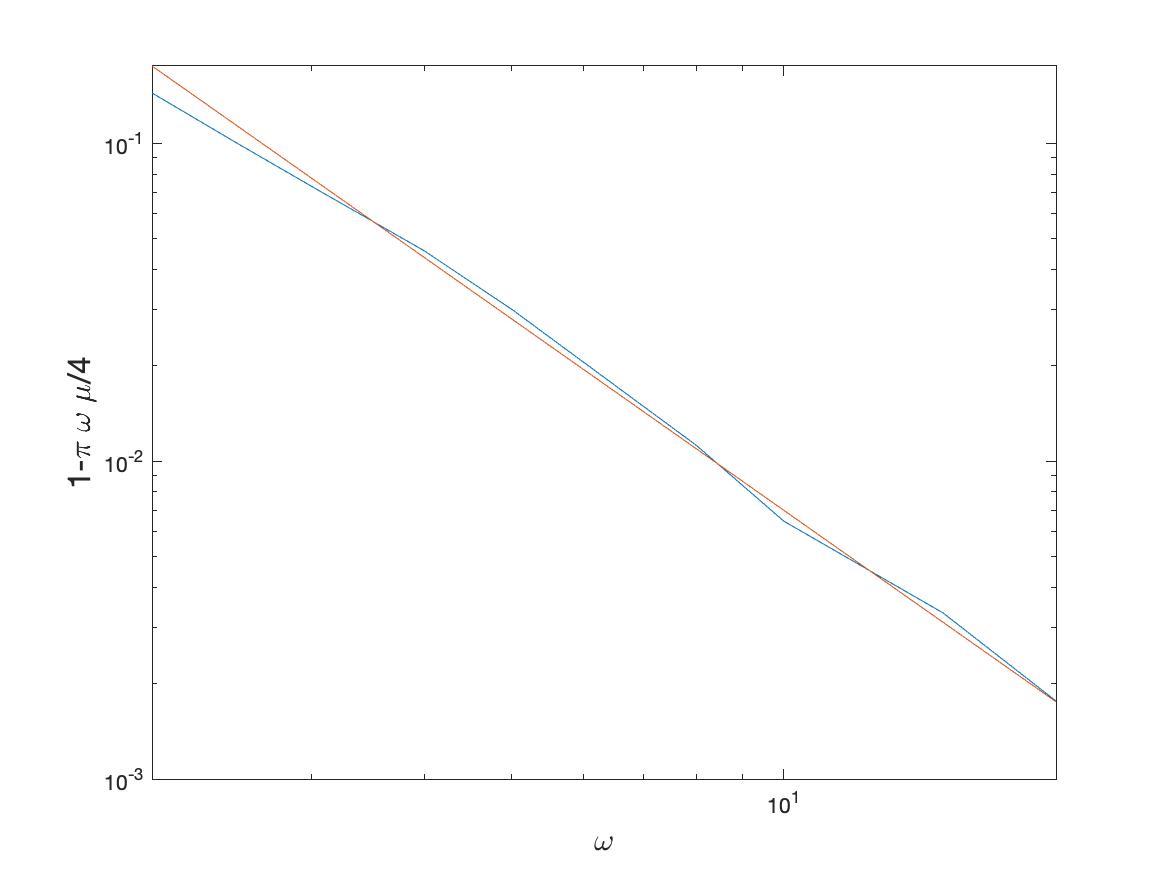}
    \caption{A log-log plot of $1 - \pi \omega \mu_{CF}/4A$ as a function of $\omega$, compared to $0.7/\omega^2$  }
    \label{fig:dec20}
\end{figure}

\noindent {\em General $\omega$:} The numerical calculations of $\mu_{CF}$ for more general $\omega \in [0,10]$ are plotted in Figure \ref{fig:4b}, and we compare the numerically computed value of $\mu_{CF}$ with the two asymptotic estimates (\ref{2}) and (\ref{2.2}), taking $L = 0.7$. Also plotted is $\mu_G$. We see that the agreement between the numerical and asymptotic calculations is excellent, with $\mu_{CF}$ lying very close to the maximum of the two asymptotic estimates for large and small $\omega$. 
\begin{figure}[htb!]
    \centering
    \includegraphics[width=0.6\textwidth]{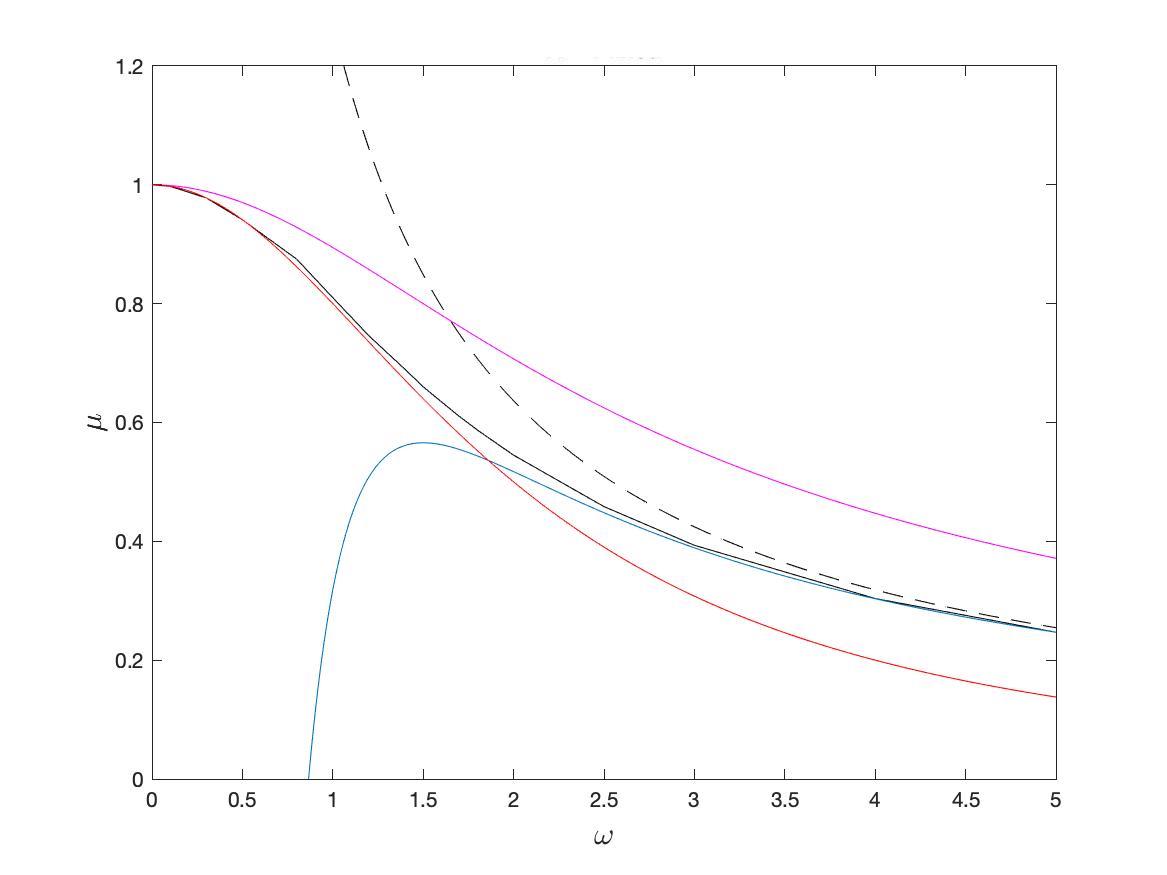}
    \caption{The location of the cyclic fold points plotted in black for $A = 1$ with $\mu_{CF}$  plotted as a function of $\omega$. In red is the small $\omega$ estimate $\mu_{CF}= 1/(1 + \omega^2/4)$ and in blue the large $\omega$ estimate $\mu_{CF} = (4A/(\pi \omega)) ( 1 - 0.7/\omega^2)$. Also plotted (dashed) is the 'averaging' estimate $\mu_{CF} \sim 4A/(\pi \omega)$ and (maroon) the grazing curve
    $\mu_G = A/\sqrt{1 + \omega^2/4}$. }
    \label{fig:4b}
\end{figure}

\vspace{0.1in}




\noindent In Figure \ref{fig:D} 
 we plot the scaled curves $(\omega \mu, \omega \langle x \rangle)$ for $\omega = 5,10,20$. The analysis presented in Appendix A implies that for larger values of $\omega$ these curves nearly coincide with $ \omega \mu_{CF} \approx 4A/\pi$. The structure of the cyclic fold is clear from this figure, with the curvature of the curves at the cyclic fold points essentially independent of $\omega$ for large $\omega$. We note from Lemma 3 that the scaled mean value at $\mu_{CF}$ given by $\omega \langle x \rangle = {\cal O}(1/\omega),$ and hence, as can be seen, takes smaller values for the larger values of $\omega$.
 


\begin{figure}[htb!]
\centering
    \includegraphics[width=0.6\textwidth]{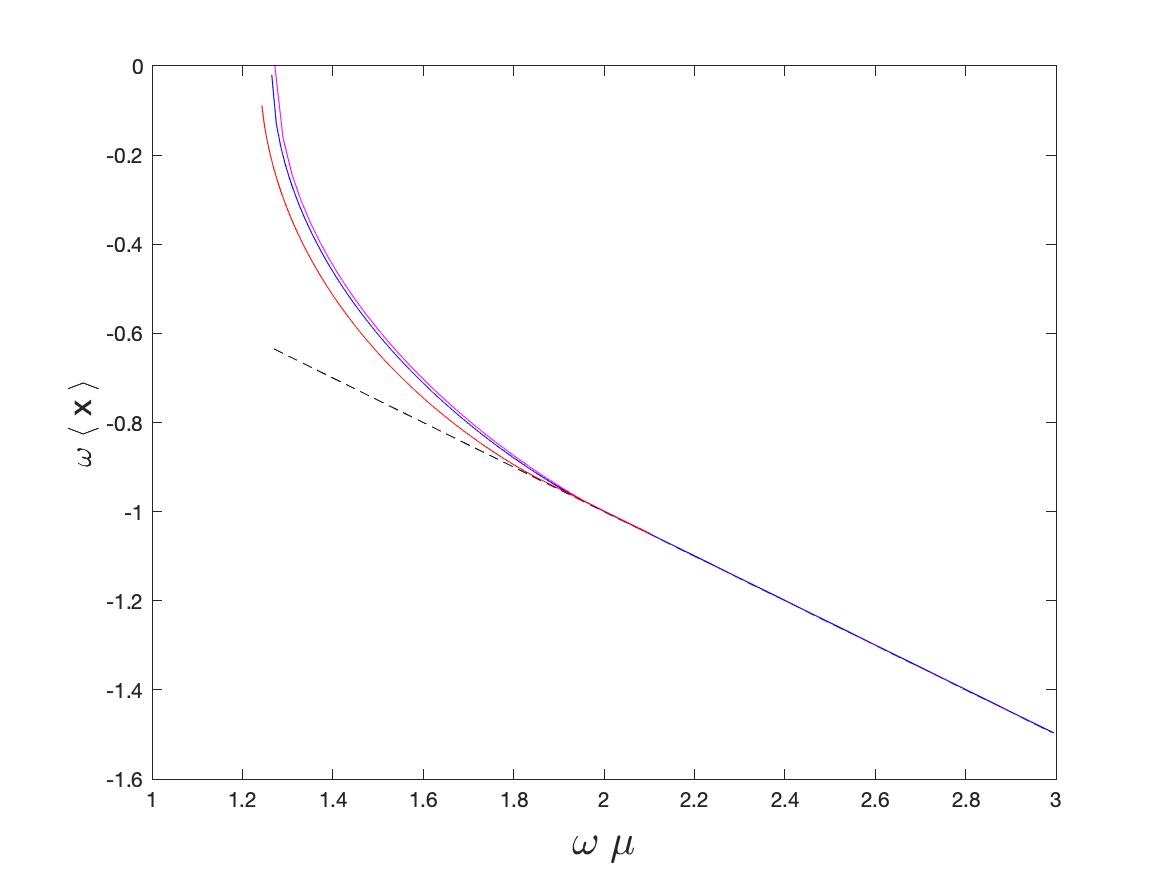}
    \caption{The value of $\omega \langle x \rangle$ plotted as a function of $\omega \mu$  with $A = 1$ and $\omega =5,10,20$ from left to right in red, blue and maroon. The line $-\omega \mu/2$ is plotted as dashed. }
    \label{fig:D}
\end{figure}

\section{Tipping under oscillatory forcing when $\mu$ has a slow drift}

\noindent In this Section we provide analytical and numerical results for tipping in the case of combined oscillatory forcing and slow drift with $d\mu/dt = -\epsilon,\  x(0) = x_0,\  \mu(0)=\mu_0$. The form of the tipping when $\epsilon = 0.1, A=1, K=10$ is shown in Figure \ref{fig:tipn} for various values of $\omega$.
\begin{figure}[htb!]
    \centering
    \includegraphics[width=0.5\textwidth]{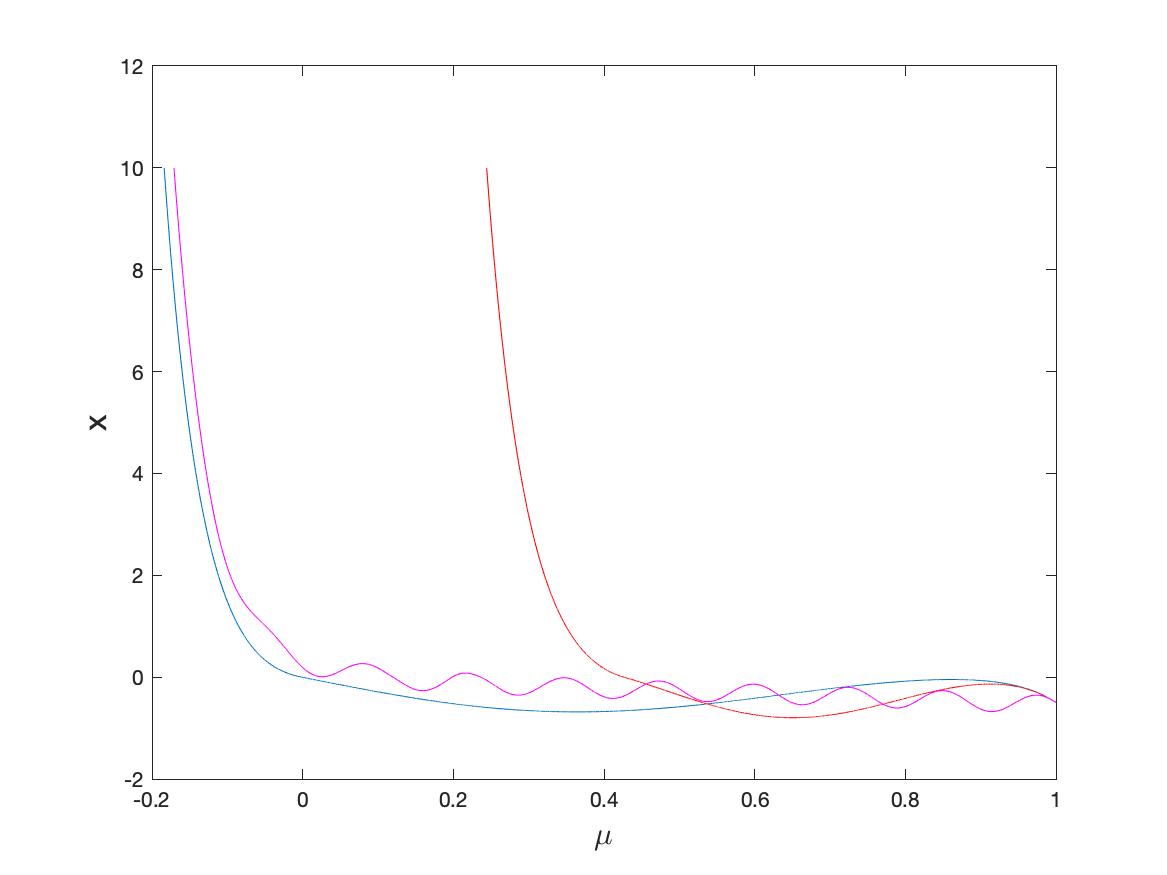}
    \caption{Tipping when $\mu_0 = 1, x_0=-1/2, \epsilon = 0.1, A = 1, K = 10$ and $\omega = 0.5, 1, 5$ (respectively in blue, red and maroon). Note that the tipping value of $\mu \equiv \mu_{TP}$ is not a monotonic function of $\omega$.} 
    \label{fig:tipn}
\end{figure}
In this section we compare and contrast the behavior of the location of the tipping point $\mu_{TP}$ for small and large $\omega$ and look at the influence of the initial value $(x_0,\mu_0)$ for smaller values of $\omega$. We also analyze the non-monotonic behavior of $\mu_{TP}$ for
smaller values of $\omega$ and make some general comments about the behavior in the $\omega - \epsilon$ plane as sketched in Figure \ref{fig:hhhhh}.

\subsection{Tipping when  $\omega$ is large.}

\noindent This case has been analysed in \cite{cody}. If the drift rate is $\epsilon$ so that
$$\frac{d\mu}{dt} = -\epsilon, $$
and $0 < \epsilon \ll 1$, a result from
\cite{cody} based on averaging for 
large $\omega$ is
\begin{equation}
\mu_{TP} \sim \mu_{CF}   - M \left( \frac{A \epsilon^2}{\omega} \right)^{1/3} \sim \frac{4A}{\pi \omega}  - M \left( \frac{A \epsilon^2}{\omega} \right)^{1/3},
\label{3}
\end{equation}
where the value of $M$ is estimated to be  $M = c_0 ( \pi/2)^{1/3}$, for $c_0$ from \eqref{1.7}. We observe that the lag in the tipping of $O(\epsilon^{2/3})$ due to slow variation through a parabolic shaped bifurcation at the cyclic fold $\mu_{CF}$ is behaving in a similar manner to tipping from a {\em static fold} (a saddle-node). Indeed, while $\mu_{CF}$ advances $\mu_{TP}$, the second term, obtained by  the canonical rescaling by $A$ and a local scaling in terms of $\epsilon^{2/3}$ \cite{haberman79} employed for the SNB, yields a lag in  $\mu_{TP}$.
\noindent We also observe that the large $\omega$ estimate (\ref{3}) for $\mu_{TP}$ is a  monotone, and smooth, function of  both $\omega$ and $\epsilon$. Furthermore it is clear that in this case $\mu_{TP} \to \mu_{CF}$ from below as $\epsilon \to 0.$


\subsection{Tipping, and transitions, for general values of $\omega$ and $\epsilon$.}

\noindent {\em Varying $\omega$ with $\epsilon$ fixed:} In the case of smaller $\omega$ we see a different behaviour in tipping with a strongly non-monotone form of the curve of tipping points as a function of $\omega$. The curves $\mu_{TP}(\omega)$ for different $\epsilon$ and $\mu_0$ are shown in the left panel of Figure \ref{fig:gg} for a set of increasing values of $\epsilon$.
\begin{figure}[htb!]
    \centering
    \includegraphics[width=0.45\textwidth]{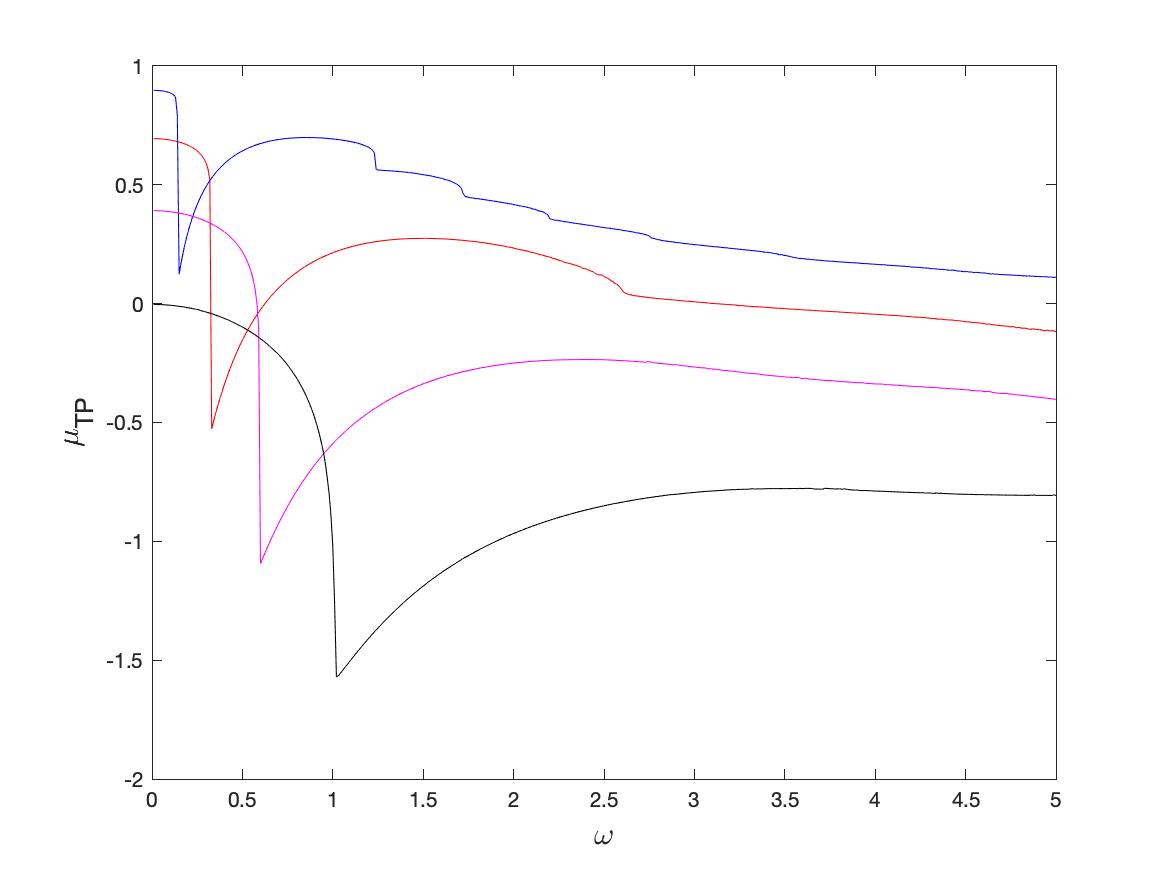}
    \includegraphics[width=0.45\textwidth]{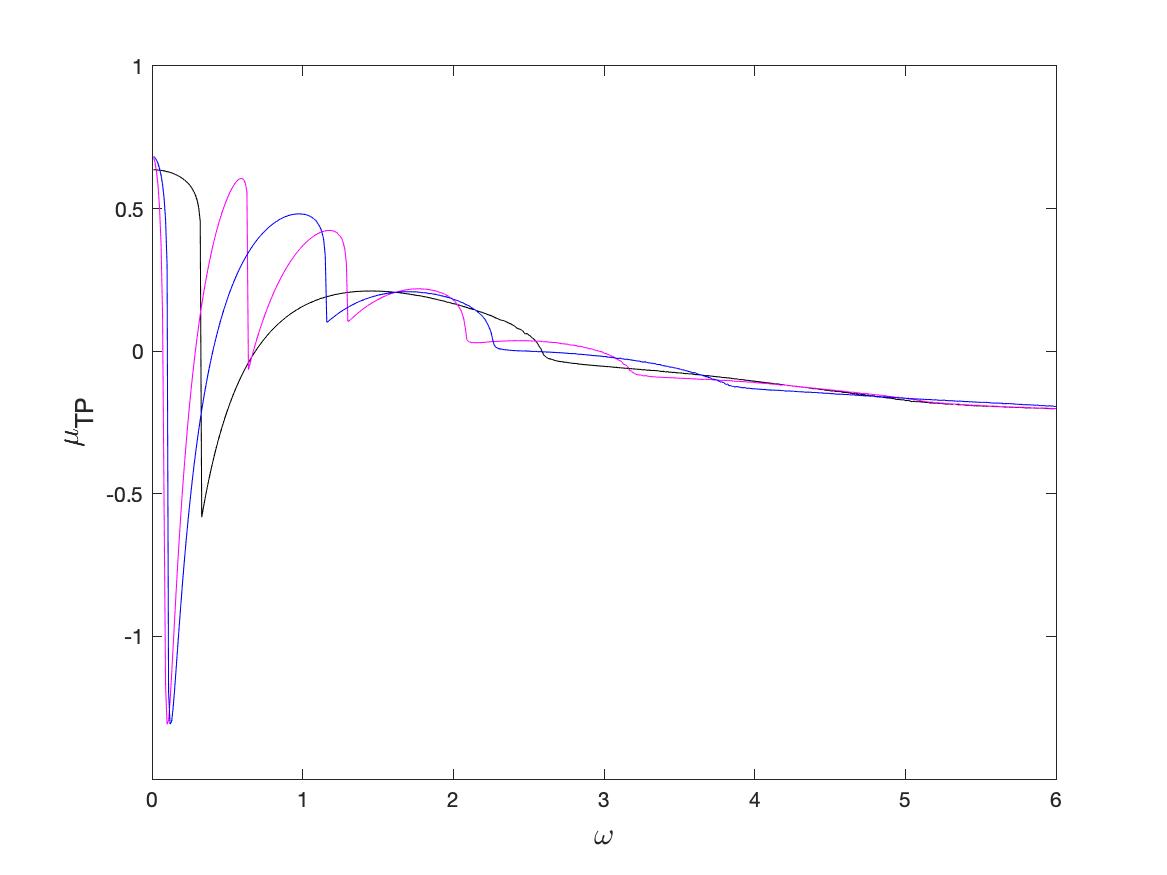}
    \caption{ Tipping points $\mu_{TP}$ when $A = 1, K = 10,  \omega \in [0,5]$, $\mu_0 = \mu_G, x_0 = -\mu_0/2$. (Left) drift rate $\epsilon = 0.025,0.1,0.25,0.5$ (respectively blue, red, maroon, black). There is a transition point at $\omega = \omega_T$ in each case, the location of which increases as $\epsilon$ increases.  (Right)  $\epsilon = 0.1$ and $\mu_0 = m \; \mu_G, x_0 = -\mu_0/2$ with $m = 1,1.5,2.$ (in black, blue and maroon).   For larger $m$ we see more transition points. }
    \label{fig:gg}
\end{figure}
These curves are revealing. For larger values of $\omega$ we see, as predicted above, smooth behaviour for all values of $\epsilon$. In contrast, for smaller $\omega$ there is a sharp transition at a value $\omega_T$, so that if $\omega < \omega_T$, near $\omega_T$ we have $\partial \mu_{TP}/\partial \omega < 0$ and large in magnitude, whilst if $\omega > \omega_T$ and close to $\omega_T$ then $\partial \mu_{TP}/\partial \omega > 0$ and is smaller in magnitude than before the transition. The value of $\omega_T$ appears to be increasing as $\epsilon$ increases. On the right panel of Figure \ref{fig:gg} we see the impact of varying the starting values $\mu_0, x_0$, by fixing $\epsilon = 0.1$ and letting $\mu_0 = m \mu_G$. For larger values of $m$ we see more transition points $\omega_T$.

\vspace{0.1in}

\noindent {\em Varying $\epsilon$ with $\omega$ fixed}. If instead we fix $\omega$ and consider $\mu_{TP}$ to be a function of $\epsilon$ then a similar pattern emerges. In Figure \ref{fig:hh}  we take $\mu_0 = \mu_G$, $\omega = 0.5$ and $\omega = 5$, and plot $\mu_{TP}(\epsilon).$ We also plot the point $(0,\mu_{CF})$ in each case. We see a number of features in this plot. Firstly that $\mu_{TP} < \mu_{CF}$, secondly, that $\mu_{TP} \to \mu_{CF}$ as $\epsilon \to 0$. We also see that for $\omega = 0.5$ then there is a sharp transition at $\epsilon \approx 0.2$, while for $\omega = 5$  there is no evidence of  such a transition.
\begin{figure}[htb!]
    \centering
    \includegraphics[width=0.5\textwidth]{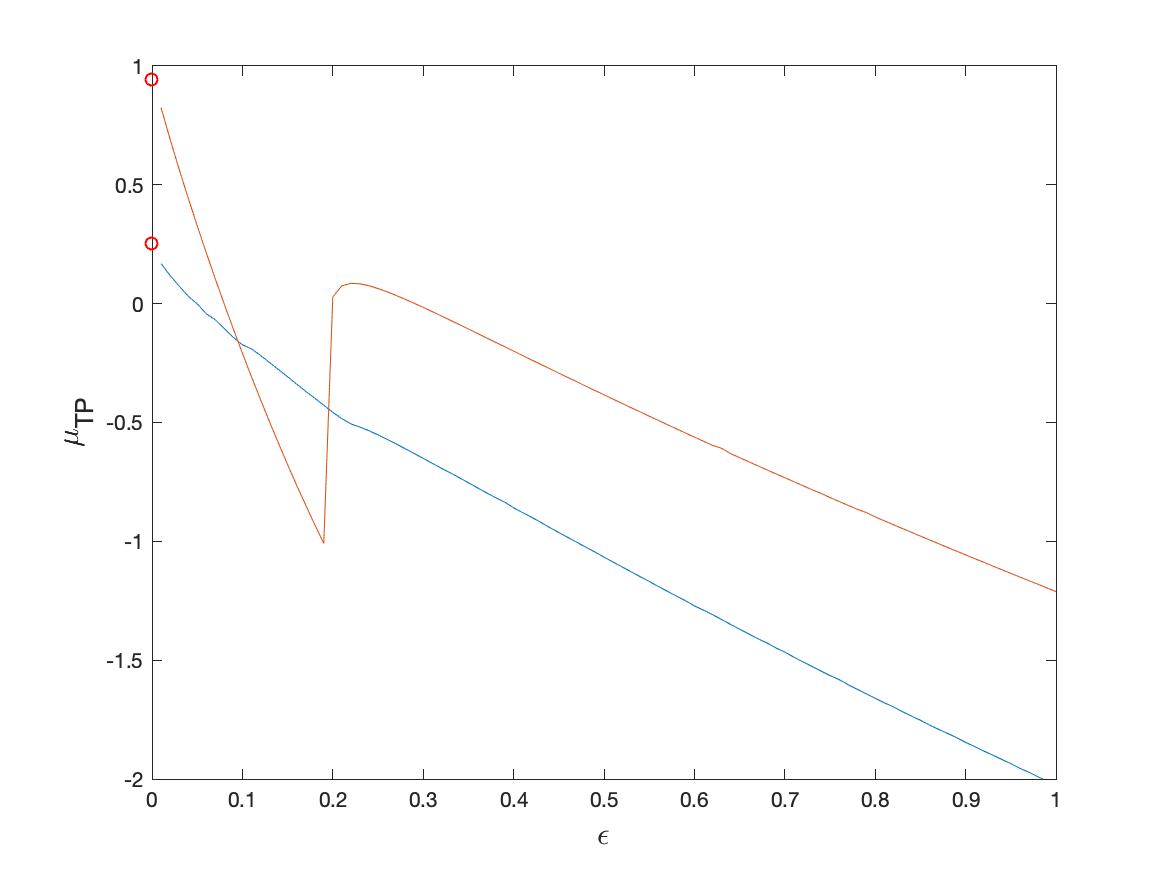}
    \caption{Tipping point $\mu_{TP}$ as a function of $\epsilon$ for $\omega = 0.5,5$ in red and blue, respectively. Also shown is the point $(0,\mu_{CF}$) (red circles). Note that (i) $\mu_{TP} < \mu_{CF}$ and (ii) $\mu_{TP} \to \mu_{CF}$ as $\epsilon \to 0$. Note that when $\omega = 0.5$ there is a large transition when $\epsilon \approx 0.2$.}
    \label{fig:hh}
\end{figure}

\vspace{0.1in}

\noindent It is clear that $\mu_{TP}$ is not a monotone decreasing function of $\epsilon$. However, we see from this figure that  $\mu_{TP} \to \mu_{CF}$ as $\epsilon \to 0$ with $\mu_{TP} < \mu_{CF}.$ 

\vspace{0.1in}

\noindent To systematise these observations we consider the surface given by the values of $\mu_{TP}$ expressed as a function of $(\epsilon,\omega)$, taking $\mu_0 = \mu_G$ and $x(0) = -\mu_0/2$. For large $\omega$ this surface is given by the asymptotic expression (\ref{3}) and is smooth. However, as we can see from the above figures, the surface has  a 'crease' for smaller values of $\omega = \omega_T$. This crease is associated with a sharp change in the gradient of the curve $\mu_{TP}$ as a function of either $\omega$ or $\epsilon$. The point $(\epsilon,\omega,\mu_{TP}) \equiv (\epsilon_T,\omega_T,\mu_{TP})$ where this transition occurs lies on the curve ${\cal T}$. We find from numerical experiments, that this curve is smooth, but terminates at the codimension-2 point when $(\epsilon_T,\omega_T) = (\epsilon^*,\omega^*)$ where the values of $\epsilon^*, \omega^*$  depend on $A$ and  other parameters such as the starting values. For $\epsilon > \epsilon^*$ or for $\omega > \omega^*$ we do not see a sharp transition. To the right of ${\cal T}$ the solution $(t,x(t))$ oscillates at most once before tipping. To the left of ${\cal T}$ the solution may oscillate many times before tipping (and indeed an infinite number of times as we approach the cyclic fold when $\epsilon = 0.$). This is illustrated in Figure \ref{fig:hhhhh}.  For larger values of $\mu_0$ we have a more complex surface, with more sharp transitions, and more lines analogous to the line ${\cal T}$. However, in this case we still see a smoother surface for large $\omega$ or $\epsilon$.

\begin{figure}[htb!]
    \centering
    \includegraphics[width=0.7\textwidth]{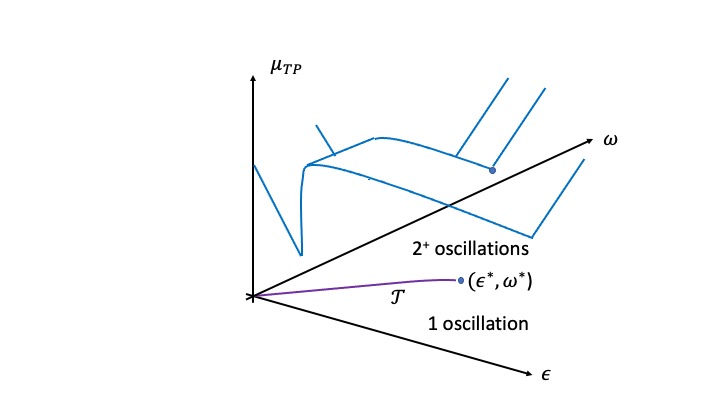}
    \caption{The surface $\mu_{TP}$ as a function of $\epsilon$ and $\omega$ when $\mu = \mu_G-\epsilon t$, $  x(0) = -\mu_G/2$. Plotted in solid is the curve ${\cal T}$ of transition points $(\epsilon_T, \omega_T, \mu_{TP})$ at which we have a transition. To the right of this curve the solution oscillates at most once before tipping and the surface drops steeply down to ${\cal T}$. To the left it can oscillate many times before tipping and the surface drops slowly down to ${\cal T}$. For $\epsilon > \epsilon^*$ or for $\omega > \omega^*$  we do not see any sharp transitions and the surface is smooth. }
    \label{fig:hhhhh}
\end{figure}

\subsection{Analysis of the behaviour of the tipping points.}\label{analysistip}

\noindent \noindent We now investigate the behaviour of $\mu_{TP}$ as a function of $\omega$ and $\epsilon$  further, and in particular study the form and reason for the sharp transitions and the termination of the curve ${\cal T}$.  We establish a series of results which support the numerical results and conjectures above.

\vspace{0.1in}

\noindent {\em The cases $\epsilon = 0$ and $\omega = 0$:}  We have already seen that $\mu_{TP} \to \mu_{CF}$ as $\epsilon \to 0$. In contrast, in the case of $\omega = 0$ then we have the ODE
$$\frac{dx}{dt} = 2|x| - \mu(t) + A.$$
This is identical to the problem with slow drift studied in Section 3, but with a simple shift in $\mu$. We therefore have
\begin{eqnarray}
    \mu_{TP} = A + {\mu}_{\epsilon}
    \label{omega0TP}
\end{eqnarray}
where ${\mu}_{\epsilon}$ is the value  estimated in the formula (2.1).

\vspace{0.2in}

\noindent {\em Increasing $\omega$ and $\epsilon$:} Consider the case of taking $x(0) = x_0 < 0$ and $\mu = \mu_0 - \epsilon t$. The solution is initially in the region $x < 0$ and can be explicitly computed. It takes the form
$$x^{-}(t) = C^-e^{-2t} - \frac{\mu_0}{2} + \frac{\epsilon t}{2} - \frac{\epsilon}{4} + \frac{A}{4 + \omega^2} \left[ \omega \sin(\omega t) + 2\cos(\omega(t)) \right], $$
for some constant $C^-$. This solution is evidently oscillatory, has an increasing mean as $t$ increases, and is exponentially stable to small perturbations. As $t$ increases then at some time $t^*$ the solution $x(t)$ crosses into the region $x > 0$, following a period where   $x = x^{-}(t)$} may have several oscillations in the region $x < 0$. 
In the region $x > 0$ we have similarly 
\begin{equation}
x^+(t) = x_P(t) + C^+ e^{2t}
\label{atl1}
\end{equation}
where $x_P(t)$ is the particular solution of the ODE when $x > 0$ given by:
\begin{equation}
x_P(t) =
 \frac{\mu(t)}{2} - \frac{\epsilon}{4} + \frac{A}{4 + \omega^2} \left[ \omega \sin(\omega t) - 2\cos(\omega(t)) \right] \qquad \mu(t) \equiv \mu_0 - \epsilon t \, .
\label{partic}
\end{equation}
The expression (\ref{atl1}) is  valid only if $x(t) > 0$ and describes a solution which is exponentially unstable.  If $x_P(t^*) < 0$ then $C^+ > 0$ and the solution rapidly increases away from $x_P(t)$ and will tip early. In contrast if $x_P(t^*) > 0$ then $C^+ < 0$ and the solution decreases away from $x_P(t)$, re-entering the region $x < 0$.  Then at a significantly later time it re-enters the region $x>0$ with tipping  following this later time and yielding a smaller value of $\mu_{TP}$ .
The value of $t^*$ depends upon $\omega$ and $\mu_0$. The numerical evidence is that $C^+ > 0$ for $\omega < \omega_T$ and $C^+ > 0$ for $\omega > \omega_T$. This shift in the tipping value underlies the 
dramatic non-monotonic behavior shown in Figure \ref{fig:gg}.
The lowest value for $\mu_{TP}$ corresponds to  taking the particular solution $x_P(t)$, extending it into the region $x < 0$, and finding the corresponding tipping point for this extended function.  

\vspace{0.1in}

\begin{figure}[htb!]
    \centering
  \includegraphics[width=0.45\textwidth]{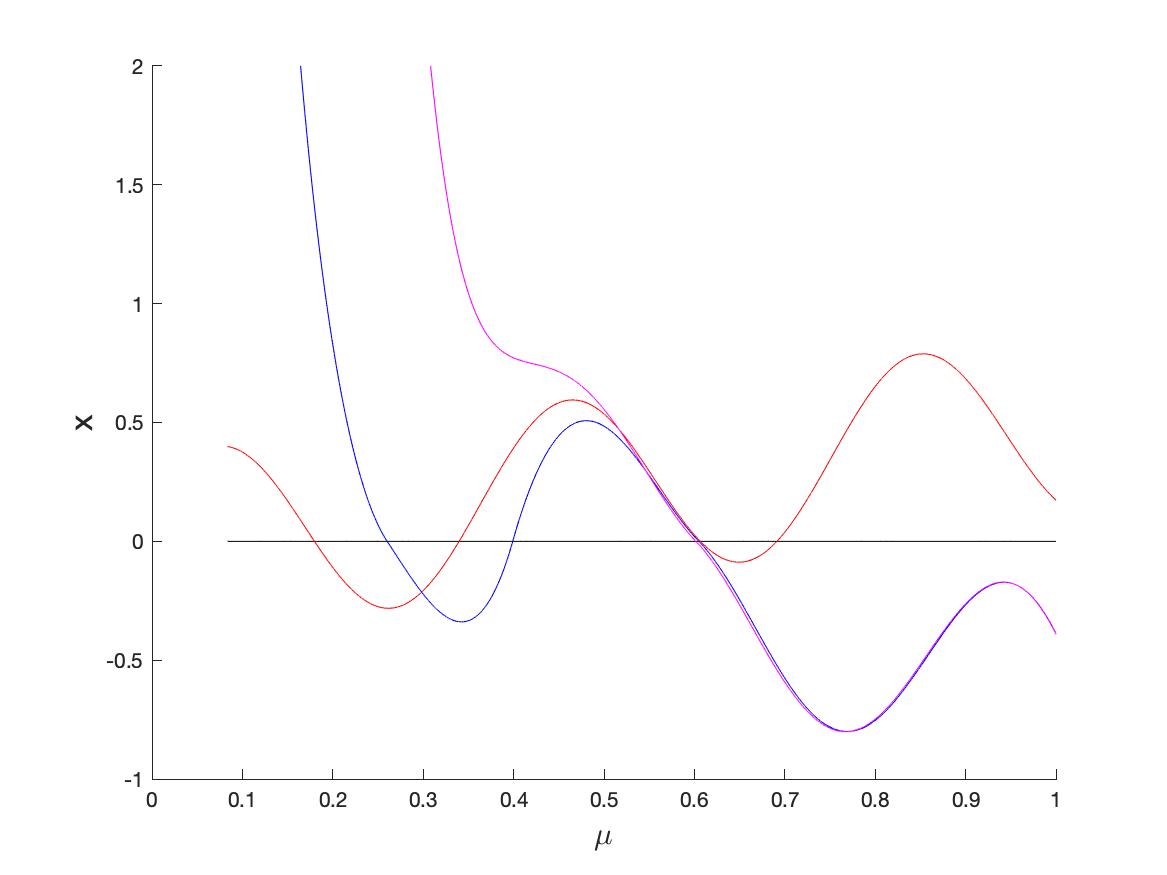}
    \caption{The solution $(\mu,x)$ when   $\omega = 1.6$ (maroon) and $\omega = 1.62$ (blue) together with the function $x_P(t)$ in red. 
    }\label{fig:22e}
\end{figure}

\vspace{0.1in}

\noindent This behaviour can be seen clearly in Figure \ref{fig:22e} where we consider two solutions of $x$ plotted as a function of $\mu$ with  $\omega = 1.6$ (maroon) and $\omega = 1.62$ (blue) together with the function $x_P(t)$ in red.  For this figure we take $\mu_0 = 1, x_0 = -\mu_0/2$.
We can see from this figure that the lower value of $\omega$ corresponds to $C^+ > 0$ with  tipping around $\mu_{TP}= 0.6$, and the larger value to tipping at the much lower value of $\mu_{TP} = 0.2.$



\noindent {\bf NOTE} Whilst the ODE (\ref{c10}) is non-smooth, the function $x(t)$ depends continuously upon the parameters $\omega$ and $\epsilon$. It follows that the tipping time $t_{TP}$ (and hence $\mu_{TP}$) must also depend continuously on these parameters. The sharp transitions in $\mu_{TP}$ as the parameters vary are thus evidence of high gradients, rather than discontinuities, in the curves $\mu_{TP}(\omega)$. Analysis of these transitions is given in  Section \ref{phase}, where we also discuss analogous transitions for $\mu_{TP}(\omega)$ also seen in the SNB context for smaller values of $\omega$.

\vspace{0.2in}

\vspace{0.1in}

\noindent Conversely, if $x_P(t) < 0$ for all $t$ then there are no sharp transitions in $\mu_{TP}$. Indeed 
$$x_P(t) = \frac{\mu_0}{2} - \frac{\epsilon t}{2} - \frac{\epsilon}{4} + \frac{A}{4 + \omega^2} \left[ \omega \sin(\omega t) - 2\cos(\omega(t)) \right] < \frac{\mu_0}{2} + \frac{A}{2\sqrt{1+\omega^2/4}} - \frac{\epsilon t}{2} - \frac{\epsilon}{4}.   $$
From this we deduce the following relationship between
$\epsilon, \omega$ and $\mu_0$:

\vspace{0.1in}

\noindent {\bf Lemma 5} {\em If 
$$\epsilon > 2 \mu_0 + \frac{2A}{\sqrt{1 + \omega^2/4}}$$
then the surface $\mu_{TP}$ will not have sharp transitions. }

\vspace{0.1in}

\noindent {\em Proof} The condition on $\epsilon$ ensures that $x_P(t) < 0$ if $t > 0.$

\qed

\noindent This lemma explains in part the form of Figure \ref{fig:hhhhh} and the termination of the curve ${\cal T}$ at the point $(\epsilon^*,\omega^*).$

\subsection{The influence of the phase for tipping at a NSF and at a SNB} 
\label{phase}

 The potential for multiple sharp decreases in $\mu_{TP}$ as a function of $\omega$ is not specific to the non-smooth context of \eqref{c10}. In \cite{zhu2015tipping}, sharp gradients in the tipping point as a function of forcing amplitude were obtained for the canonical SNB with a slowly varying bifurcation parameter, and low frequency forcing, comparable to $\omega < 1$.
 Here we analyze the role of the initial condition $\mu_0$, which can be interpreted as phase, together with $\epsilon$ in determining the size of the fluctuations and the number of jumps seen in $\mu_{TP}$ as a function of $\omega$.  We discuss how the method can be applied also in the SNB context. The results speak to challenges in predictability of tipping, varying considerably with the initial state, in the setting of $\omega \ll 1$ for both the smooth and the non-smooth contexts.

\subsubsection{ The role of $x_P(t)$ and forcing }

As shown in Section \ref{analysistip}, the relative crossings of $\Sigma$ 
 by $x(t)$ and the particular solution $x_P(t)$ for  $x>0$ \eqref{partic}  feature centrally in the analysis of the location of the tipping points $\mu_{TP}$. Recall that
$x_P(t)$ has bounded oscillations, and decreases as $t$ increases (and hence $\mu$ decreases). Consequently it is negative for sufficiently large $t$. If $x(t) > 0$  when $x_P(t) < 0$ then the coefficient $C^+$ in \eqref{atl1} is positive, and hence $x(t)$ tips. In contrast if $x_P (t)> 0$ as $x(t)$ crosses $\Sigma$ into the region $x>0$ then we have $C^+ < 0$ and the solution tips at a later time when $ x_P(t) < 0$. This comparison of  the relative phase of $x(t)$ and $x_P(t)$ as they cross $\Sigma$ leads to analytical expressions that capture the sharp transitions  of $\mu_{TP}$ for $\omega$ not large.  While $x_{P}(t)$ is a convenient function to use for \eqref{c10}-\eqref{c11.5}, obtaining an explicit form for a particular solution does not generalize easily to nonlinear settings, e.g.  such as the polynomial form of the canonical SNB model.  Instead we give the analysis  in terms
$-\mu(t) + f(t)$ \eqref{c11.5},
which is then generalizable in nonlinear settings.

 \begin{figure}[htb!]
    \centering
    \includegraphics[width=0.45\textwidth]{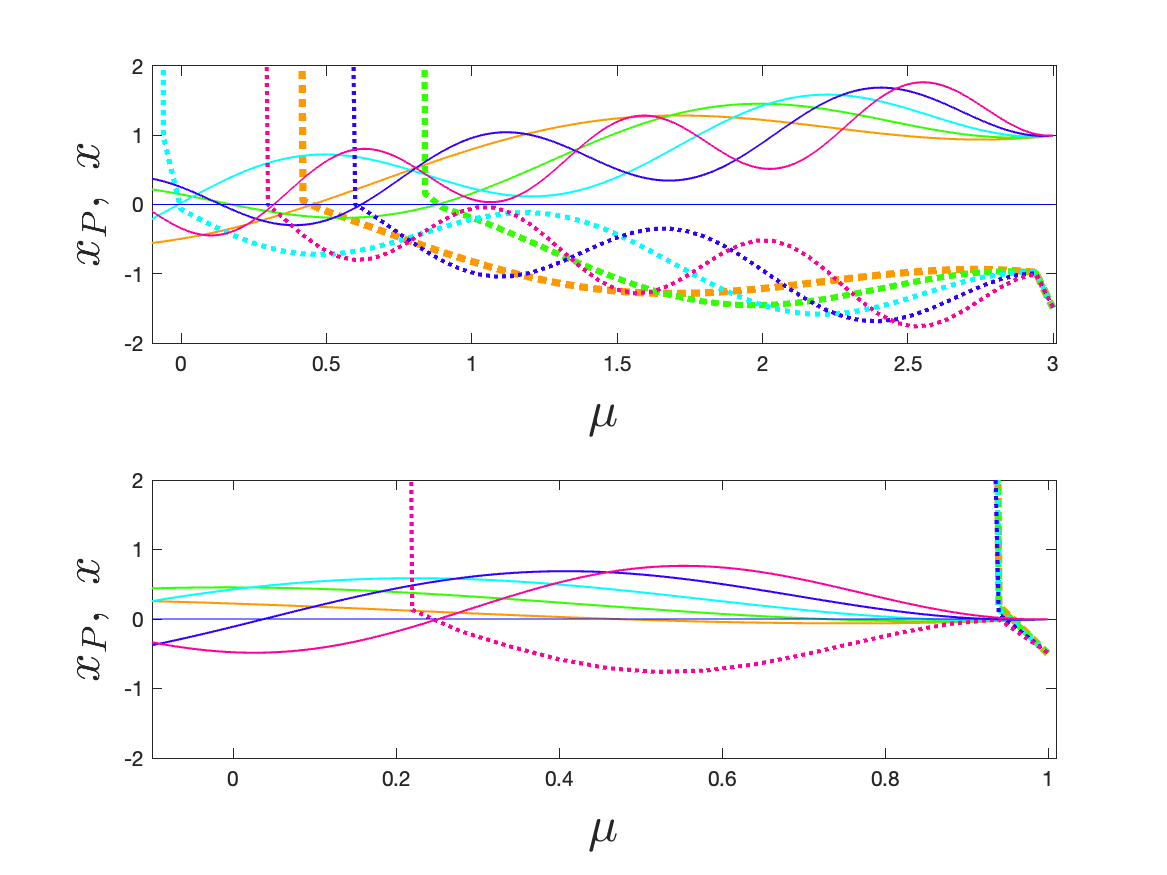}
\includegraphics[width=0.45\textwidth]{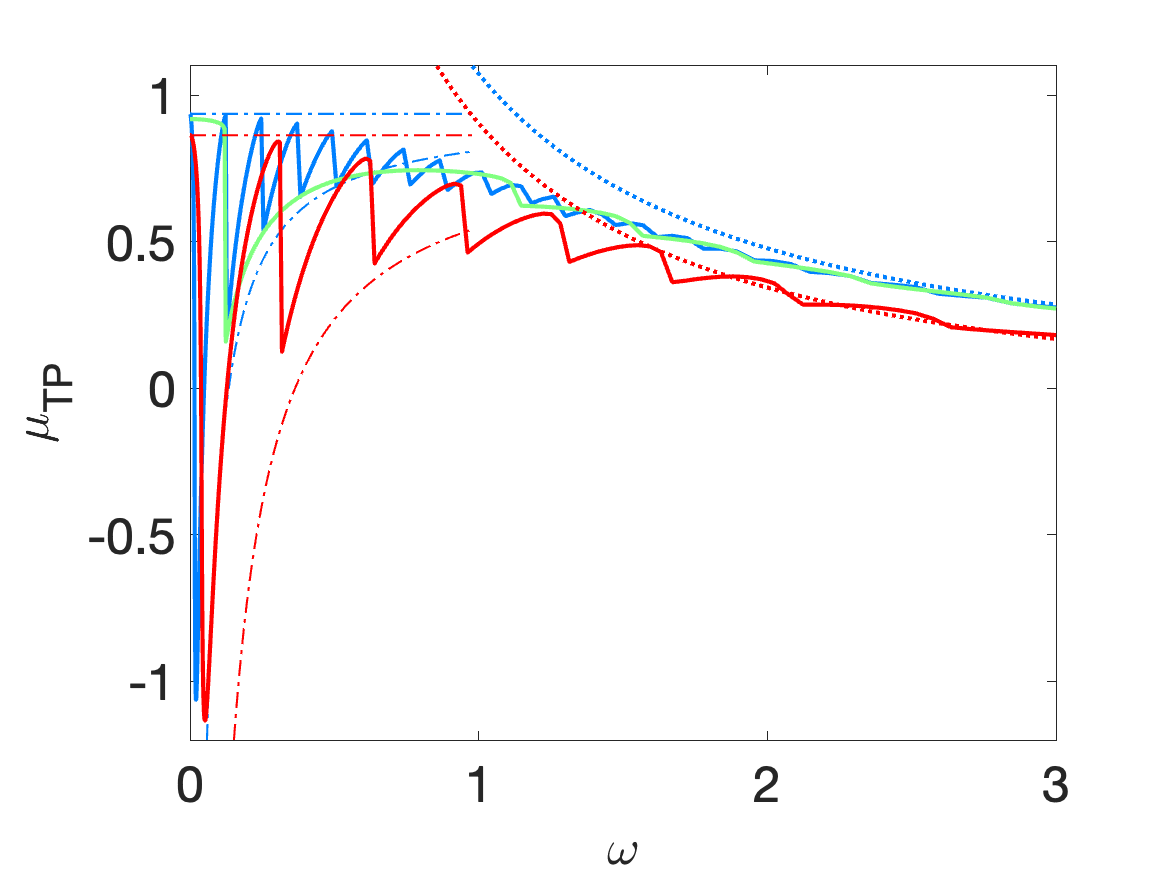}\\
    \hskip 1 in  (a) \hskip 3in (b) 
    \caption{(a) Dotted (solid) lines show $x(t)$ ($x_{P}(t)$)  as functions of $\mu$ for five values of $\omega$,  $0.04<\omega<0.13$, and $\epsilon=0.02$,  with $\mu_0 = 3\mu_G$ (upper) and $\mu_0 = \mu_G$ (lower). The values of $\omega$ increase with the order of colors orange, green, cyan, purple, pink. For $\mu_0 = \mu_G$, the first four trajectories for $x(t)$ are nearly indistinguishable, in contrast to larger $\mu_0$.
    (b) Solid lines show $\mu_{TP}$ vs. $\omega$  for
    $\epsilon = 0.02$ and $\mu_0 = \mu_G$ (green), $\epsilon = 0.02$ and $\mu_0 = 2\mu_G$ (blue), and  $\epsilon = 0.05$ and $\mu_0 = 2\mu_G$ (red). 
    Multiple sharp transitions  in $\mu_{TP}$ are observed for  $\omega <1$ for larger $\mu_0$, with a  range that decreases with $\omega$. These fluctuations continue for $\omega>1$,  approaching \eqref{3} (dotted line) for $\omega >2$.
     Note \eqref{3} is independent of $\mu_0$, so green and blue dotted lines are identical. Dash-dotted lines show upper and lower bounds on $\mu_{TP}$ in the cases with $\mu_0 = 2\mu_G$,
     as given by \eqref{muTPbound}  for $\omega<1$ .}
    \label{phasefigs}
\end{figure}
  
\vspace{0.1in}

 \noindent There are different factors contributing to the number and size of the transitions in $\mu_{TP}$, which, as one might expect from Lemma 5, are  commonly observed for small $\epsilon$.
 For example, for increasing $\mu_0$, the function $x_P(t)$ is {\em positive} over a larger interval of $t$, which affords more opportunity for a mismatch in the relative phase of $x(t)$ and $x_P(t)$ as they cross $\Sigma$. Figure \ref{phasefigs} (a) illustrates the variation in the relative phase of $x(t)$ and $x_P(t)$ in crossing $\Sigma$ over a range of frequencies $\omega<1$,  contrasting the case of taking $\mu_0 = \mu_G$, where there is less variation in $\mu_{TP}$,  as compared  to the repeated non-monotonic behavior of $\mu_{TP}$ over the same values of $\omega$ for the larger value of $\mu_0= 2 \; \mu_G$. The latter indicates multiple sharp transitions in $\mu_{TP}$
 for $\omega<1$ and larger values of $\mu_0$. Likewise for $1<\omega=O(1)$, with higher frequency oscillations, there are additional opportunities for phase mismatch yielding fluctuations in $\mu_{TP}$.
 In general,  as the period $T$ of  the function $x_P(t)$ is $T = 2\pi/\omega$, the  sharp transitions seen in $\mu_{TP}$  are of order 
 \begin{eqnarray}
\Delta \mu_{TP} \propto 2\pi \epsilon/\omega. \label{DeltamuTP}
\end{eqnarray}
Hence, the size of these transitions decreases in magnitude as $\omega$ increases.
Note that these transitions correspond to fluctuations around the analytical approximation \eqref{3} for $\omega>2$, with the range of fluctuations reducing with increasing $\omega$. 
Figure \ref{phasefigs}(b) illustrates these general characteristics of the jumps in $\mu_{TP}$, which are analyzed below.\\

\noindent {\bf Analysis of non-monotonic $\mu_{TP}$ for a NSF:}\\

\vspace{0.1in}

\noindent The analysis of \cite{zhu2015tipping} determines the tipping of the dynamic SNB solution in terms of
roots of the time-varying non-autonomous forcing in SNB, analogous to $-\mu(t) + f(t)\equiv g_{NS}(\mu)$ as in \eqref{c12},\eqref{c11.5}. Writing the roots of $x_P(t)$ and $g_{NS}(\mu)$ in terms of $\mu$ and denoting these as $\mu_P$ and $\mu_r$,
respectively, we note that $\mu_r = \mu_P + O(\epsilon)$. Then we consider the roots $\mu_r$ such that $g_{\rm NS}(\mu_r)=0$ and $g_{\rm NS}'(\mu_r)<0$, which correspond to $\mu_P$ where $x_P(\mu_P) = 0$ and $\frac{d x_P}{d\mu}(\mu_P)>0$. 
These values yield the potential for tipping as $x_P$ crosses below $\Sigma$, and thus $\mu_r$ serves as a surrogate for $\mu_{TP}$ away from its steep drops in value. Indeed, for the SNB setting in \cite{zhu2015tipping}, the tipping value is determined as $\mu_r+ {\cal O}(\epsilon^{2/3}) $. In the following we assume that $\epsilon$ (and $\omega$) are sufficiently small, so that Lemma 5 does not apply, and we expect to see sharp transitions in the location of the tipping point $\mu_{TP}$ and hence of $\mu_r$.


\noindent Following a formulation of \eqref{c10}  similar to \cite{zhu2015tipping}, we write  \eqref{c11.5} in terms of $\mu$, rather than $t$:
\begin{eqnarray}
 & &-\epsilon x_\mu = g_{\rm NS} + 2|x|\\
 & &g_{\rm NS}(\mu) = -\mu + A\cos(\omega t) = 
-\mu + A\cos(\Omega(\mu_0-\mu)) \qquad \omega = \epsilon \Omega\label{gns_def}
 \end{eqnarray}
For $g_{\rm NS}(\mu_r)=0$ and $g_{\rm NS}'(\mu_r)<0$ we have
\begin{eqnarray}
\frac{\mu_r}{A} & = & \cos(\Omega(\mu_0-\mu_r)) \implies \frac{\mu_r}{A}\leq 1 \label{eqmur}\\
1 & >  & \Omega A\sin(\Omega(\mu_0-\mu_r)) \, . \label{dgmurpos}
\end{eqnarray}
    Differentiating  \eqref{eqmur} we obtain the behavior of $\frac{\partial \mu_r}{\partial \Omega}>0$
\begin{eqnarray} 
   \frac{\partial \mu_r}{\partial \Omega}\left[ 1- \Omega A\sin(c(\mu_0-\mu_r))\right] &=& - (\mu_0-\mu_r) A\sin(\Omega(\mu_0-\mu_r)) \ \ \implies \nonumber\\
   \frac{\partial \mu_r}{\partial \Omega} &=&-\frac{ (\mu_0-\mu_r) A\sin(\Omega (\mu_0-\mu_r))} { 1- \Omega A\sin(\Omega (\mu_0-\mu_r))} \label{sgndmurdc}
   \end{eqnarray}

\noindent The sign of $\frac{\partial \mu_r}{\partial \Omega}>0$ follows from  \eqref{dgmurpos} and the corresponding phase of the oscillations in $g_{NS}$. From \eqref{sgndmurdc}, we conclude that $\mu_{r}$ increases with $\omega$, with a rate that increases with $\mu_0$, as illustrated by comparing the two cases for $\epsilon = .05$ in Figure \ref{phasefigs}(b). 
From \eqref{eqmur}, $\mu_r$ cannot increase to values greater than $A+\mu_\epsilon$ in \eqref{omega0TP}. 
Then as $\mu_0$ increases, there must be  multiple intervals of $\frac{\partial \mu_r}{\partial \Omega}>0$ in order to be consistent with \eqref{sgndmurdc}.
Figure \ref{phasefigs}(b) shows approximate upper and lower bounds 
 (dash-dotted lines) for the range of values that $\mu_r$ takes for $\omega<1$ and different combinations of $\Omega$ and $\epsilon$. These are given by
\begin{eqnarray}
    A+\mu_\epsilon< \mu_{TP} < A+\mu_\epsilon - \frac{2\pi\epsilon}{\omega}
    \label{muTPbound}
\end{eqnarray}
obtained by combining the results above \eqref{omega0TP} and \eqref{DeltamuTP}.
 
\vspace{0.1in}

\noindent To provide some additional quantitative insight into these intervals of $\frac{\partial \mu_r}{\partial \Omega}>0$ , we then rewrite \eqref{eqmur} in the form, 
\begin{eqnarray}
    \Omega = \frac{\cos^{-1}(\mu_r/A)}{\mu_0-\mu_r - \frac{2n\pi}{\Omega}}, \qquad n \in {\mathbf Z}, n \geq 0 .\label{ceqmur}
\end{eqnarray}
For small $\omega$,
 e.g. $\omega <1$ for which there are sharp transitions in the tipping point, we can relate $n$ to the number of intervals of $\frac{\partial \mu_r}{\partial \Omega}>0$ in an interval $0<\omega<\omega^*$, assuming $\mu_0> A+\mu_\epsilon$. Note that for sufficiently small values of $\epsilon$, (see Lemma 5) there is always at least one value of $\omega<1$ for which there is a sharp decrease in $\mu_r$, given \eqref{omega0TP} for $\omega=0$ and a sharp decrease of $O(\omega^{-1})$ \eqref{DeltamuTP} for small $\omega>0$. 
    Specifically, as $\mu_r$ approaches $A +\mu_\epsilon$  we  approximate the value $\omega^*$ at which there is a sharp transition by setting $\mu_r=A+\mu_\epsilon$ in \eqref{ceqmur},
    \begin{eqnarray}
      & &  \omega^*= \epsilon
      \frac{\cos^{-1}(1+\mu_\epsilon/A)} 
       {\mu_0-A-\mu_\epsilon} \sim \frac{2n\pi + O(\sqrt{\mu_\epsilon})}{\mu_0-A-\mu_\epsilon } \, . \label{murneg}
        \end{eqnarray}
   We can use this expression to determine the number of sharp transitions of $\mu_{TP}$ over a range of small $\omega$.  Then $n_{max}$, defined  as 
  \begin{eqnarray}
n_{max} & =&  \underset{n}{\arg\min} \ Q(n; \omega) \mbox{ \  where  \ }\nonumber\\
 & &Q(n;\omega) =  \omega - \epsilon
      \frac{\cos^{-1}(1+\mu_\epsilon/A)} 
       {\mu_0(\omega)-A-\mu_\epsilon}
       =  \omega - \epsilon
      \frac{2n\pi + O(\sqrt{\mu_\epsilon})} 
       {\mu_0(\omega)-A-\mu_\epsilon}>0 \, ,
      \label{nmax}
  \end{eqnarray}
  approximates the number of intervals of  $\frac{d\mu_r}{d\Omega}>0$ for frequencies below a given (small) value of $\omega$.
  Note that  $\mu_0$ may depend on $\omega$, e.g. as above where $\mu_0$ is an integer multiple of $\mu_G$.
  Illustrations of $n_{max}$  appear in Figure \ref{phasefigs}(b)    where 
  \eqref{nmax} yields $n_{max} =8$ and $n_{max}=3$ for $\omega = 1$ and $\mu_0 = 2\mu_G$
  with $\epsilon = 0.02$ and $\epsilon = 0.05$, respectively. As mentioned above, a larger number of sharp drops between intervals with large $\frac{d\mu_r}{d\Omega}$ suggests a strong sensitivity to the initial conditions when predicting future tipping.

\vspace{0.1in}
   
\noindent {\bf Tipping at a SNB}\\

\begin{figure}[!htbp]
	\centering
		  \includegraphics[width=3.in]{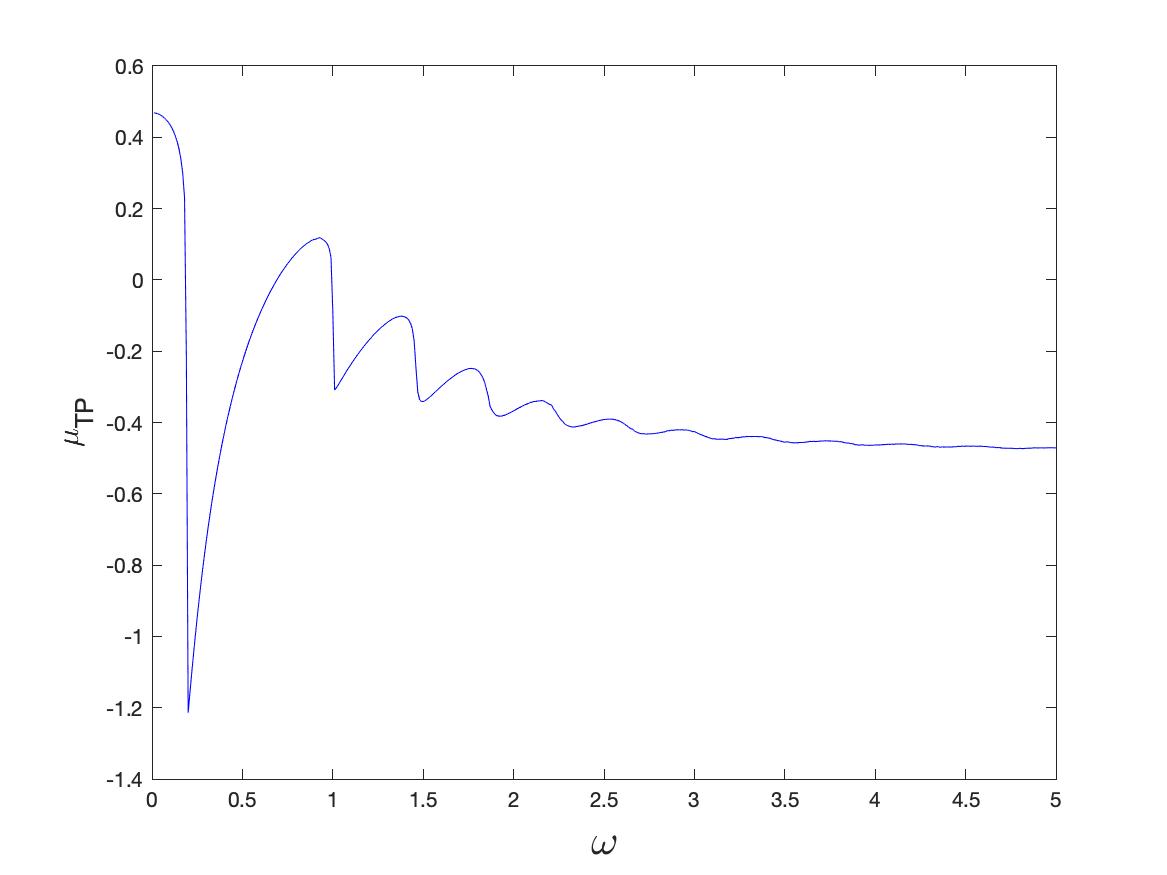}
    \includegraphics[width=3.in]{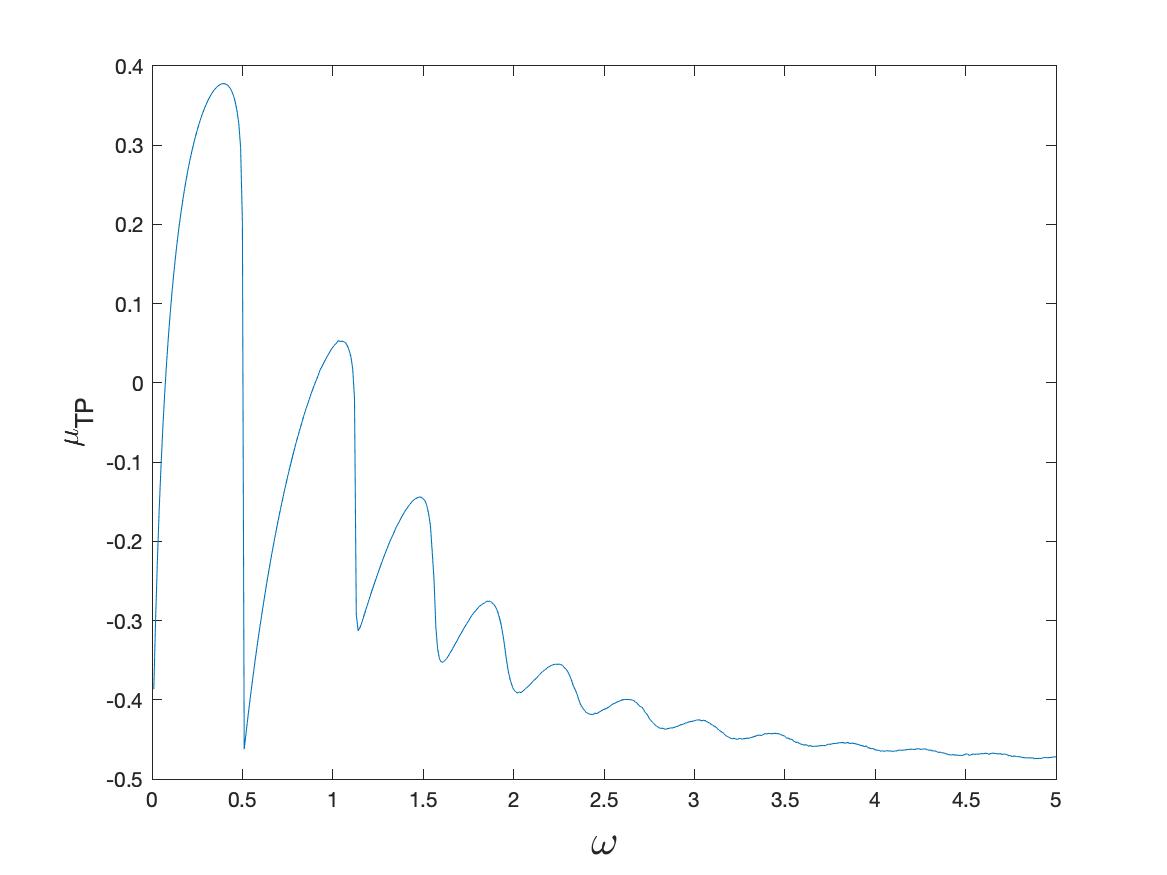}
	\caption{Calculated values of $\mu_{TP}$ for the SNB problem when $\epsilon = 0.1, K = 10$ and $A = 1$. Here we take $\mu = 1 - \epsilon t$ and $x(0) = -1$. 
 On the left we consider the problem $dx/dt = x^2 - \mu + A \cos(\omega t)$ and on the right the problem $dx/dt = x^2 - \mu + A \sin(\omega t)$. We compare this figure with Figure 1 for the NSF. See also the corresponding  Figure \ref{fig:22ff} for the smoothed NSF.} 
	\label{fig:mar1}
\end{figure}

\noindent We recall the approach from \cite{zhu2015tipping} for the canonical smooth SNB model.
For small $\omega$, it is again convenient to write the equations in terms of $\mu$ rather than $t$, defining a function  $f_{\rm SNB}$
\begin{eqnarray}
& &dx/dt = x^2 - \mu + A\cos(\omega t), \qquad
d\mu/dt = -\epsilon, \ \ \mu(0) = \mu_0 \nonumber\\
 & &-\epsilon dx/d\mu = f_{\rm SNB} + x^2\nonumber\\
 & &f_{\rm SNB}(\mu) = -\mu + A\cos(\omega t) = 
-\mu + A\cos(\Omega(\mu_0-\mu)) \label{fsnb_def}
\end{eqnarray}
An outer solution for $x(\mu)$  motivates a local expansion near $\mu_r$ such that
\begin{eqnarray} 
x(\mu_r)=0 \mbox{ \  for \  }  f_{\rm SNB}(\mu_r) = 0, \mbox{ and }
f_{\rm SNB}'(\mu_r) = O(1).
\end{eqnarray}
Similar to the analysis above for the NSF case, one can analyze the behavior of $f_{SNB}$ to obtain the sequence of sharp transitions in $\mu_{SNB}$ for small $\omega$, which we leave as an exercise.  The behavior is shown in Figure \ref{fig:mar1} for the case of the forcing $A \cos(\omega t)$ on the left, and $A \sin(\omega t)$ on the right. In both cases we take $\mu = 1 - \epsilon t$ and $x(0) = -1$ with $\epsilon = 0.1$. There are clear similarities between this figure for the SNB and the corresponding Figure \ref{fig:dec1} for the NSF.   Given the similar behavior observed for both NSF and SNB, it follows that $\mu_{TP}$ for the smoothed NSF shows the same characteristics, as shown in  Figure  \ref{fig:22ff}.


\section{The impact of smoothing}

\noindent In this section we will briefly extend the results obtained for the slow drift and the oscillatory forcing from the non-smooth system (\ref{c10}) to the smoothed system (\ref{c11}).

\subsection{General results}

For general forcing, the affect of smoothing is always to {\em postpone} tipping.

\vspace{0.2in}

\noindent {\bf Lemma 6} {\em If the smoothing value is $\alpha> 0$ and $x(t,\alpha)$ satisfies the  smoothed  equation
\begin{equation}
\frac{dx}{dt} = 2\sqrt{x^2 + \alpha^2} - 2 \alpha - \mu + A  \cos(\omega t), \quad \frac{d\mu}{dt} = -\epsilon, \quad x(0) = x_0, \  \mu(0) = \mu_0
\label{dec2}
\end{equation}
with $x_0,\mu_0$ independent of $\alpha$. Then the value of $t \equiv t_{TP}$ at tipping {\em increases} with $\alpha$, and hence  $\mu_{TP}$ at tipping {\em decreases} with $\alpha$. }

\vspace{0.1in}

\noindent {\em Proof} Differentiating (\ref{dec2})  with respect to $\alpha$ we have that $x_{\alpha} \equiv \partial x/\partial \alpha$ 
 satisfies the differential equation
$$\frac{d x_{\alpha}}{dt} = \frac{2 x x_{\alpha}}{\sqrt{x^2 + \alpha^2}} - 2, \quad x_{\alpha}(0) = 0.$$
Hence, as $x'_{\alpha}(0) = -2$ it follows that $x_{\alpha}$ is negative for $t>0$ provided that $t$ is sufficiently small. Suppose that there is some first later time $t_1$  at which $x_{\alpha}(t_1) = 0$,
then  at $t_1$   we must have $dx_{\alpha}/dt \ge 0$. 
However from the above differential equation we have $dx_{\alpha}/dt < 0$ when $x_{\alpha} = 0$. We conclude that there is no such time $t_1$ and hence that $x_{\alpha}(t) < 0$ for all $t > 0.$

\vspace{0.1in}

\noindent Now suppose that $x(t_{TP}, \alpha) = K$. Differentiating with respect to $\alpha$ we have
$$x_{\alpha}(t_{TP},\alpha) + \left.\frac{dx}{dt} \right|_{t=t_{TP}}\frac{\partial t_{TP}}{\partial \alpha} = 0.$$
Hence
$$\frac{\partial t_{TP}}{\partial \alpha} = -x_{\alpha}/\left(\left. dx/dt \right|_{t=t_{TP}}\right).$$
Now, we have that $x_{\alpha} < 0$, and also by the definition of $K$ it follows that  $\left. \frac{dx}{dt} \right|_{t=t_{TP}} > 0$. We deduce that
$$\frac{\partial t_{TP}}{\partial \alpha} > 0.$$
Now $\mu_{TP} = \mu_0 - \epsilon t_{TP}.$ Hence $t_{TP}$ increases with $\alpha$ and $\mu_{TP}$ decreases.

\qed

\subsection{Smoothed slow drift}

\noindent We first consider the case of the unforced system with slow drift $\epsilon$. The value of the tipping point $\mu_{TP}$ here depends upon the  balance between the smoothing parameter $\alpha$ and the drift rate  $\epsilon$. If the drift rate is greater than the smoothing, then the system is dominated by the non-smooth behavior. Conversely if the system drifts slowly then smoothing becomes more important. We also see that the smooth saddle-node estimate (\ref{1.7}), is a good approximation  only for a range of $\alpha$ bounded below by $\epsilon$ and above by a function of $K$ and $\epsilon$.

\vspace{0.1in}

\noindent We first consider the limits of the tipping behaviour as $\alpha \to 0$, and as $\alpha \to \infty$.

\vspace{0.2in}

\noindent {\bf Lemma 7} If $\mu = 1 - \epsilon t, x(0) = 0$ then

\vspace{0.1in}

\noindent {\em (i) As $\alpha \to 0 $ 
\begin{equation}
\mu_{TP}(\alpha) \to - \epsilon \log(2K/\epsilon)/2 \equiv \mu_{TP}(0).
\label{apr1}
\end{equation}

\vspace{0.1in}

\noindent (ii) As $\alpha \to \infty$ 
\begin{equation}
\mu_{TP}(\alpha) \to -\sqrt{1 + 2 \epsilon K} \equiv \mu_{TP}(\infty).
\label{c41}
\end{equation}

\vspace{0.1in}

\noindent (iii) $\mu_{TP}(0) > \mu_{TP}(\alpha) > \mu_{TP}(\infty). $
}

\vspace{0.2in}

\noindent {\em Proof} The proof of (i) follows immediately from continuity arguments, and from the value of $\mu_{\epsilon}$ in \eqref{ca7}.

\vspace{0.1in}

\noindent To prove (ii) we note that if $|x| < K$ then as $\alpha \to \infty$,  $2\sqrt{\alpha^2 + x^2} - 2\alpha \to 0.$ Hence, to leading order in $\alpha^{-1}$ for large $\alpha$, $x$ satisfies the simple ordinary differential equation
$$dx/dt = -\mu(t).$$
If we set $x(0) = 0$ and $\mu(t) = 1 - \epsilon t$ then we have 
$x(t) \sim \epsilon t^2/2 - t$.
Hence, $x = K$ if $t = (1 + \sqrt{1 + 2 \epsilon K})/\epsilon$,  which yields (ii).
 \vspace{0.1in}
 
 \noindent The result (iii) follows from the fact that $\mu_{TP}(\alpha)$ is a monotone decreasing function of $\alpha$.
 
 \qed

\vspace{0.1in}


\noindent We now consider the behaviour of $\mu_{TP}$ for more general values of the smoothing parameter $\alpha$. We identify three ranges for $\alpha$ that correspond to the asymptotic behaviour described in Lemma 7: (i) small $\alpha < \alpha_0$ where the 'non-smooth' limit (\ref{apr1}) applies, 
(ii) large $\alpha > \alpha_1$ where the limiting value of $x(t) = K$ is important and estimate (\ref{c41}) applies, and
(iii) intermediate $\alpha_0 < \alpha < \alpha_1$ where the usual SNB estimate \eqref{1.7} applies. 

\vspace{0.1in}

\noindent  For case (i) we consider the limit of $\alpha \to 0$. We note, that comparing the SNB estimate (\ref{1.7}) with the non-smooth estimate (\ref{apr1}) then
$$c_0 \; \alpha^{1/3} \epsilon^{2/3} < \epsilon \log(2K/\epsilon)/2, \quad \mbox{when} \quad \alpha < \epsilon (\log(2K/\epsilon)/2)^{3} \equiv \alpha_0,$$
noting that $c_0 > 1$.
This implies that the non-smooth estimate for the tipping value is more accurate if $\alpha < \alpha_0.$ 

\vspace{0.1in}

\noindent To obtain the results for cases (ii) and (iii) we note that as $\alpha \to \infty$
$$\quad 2 \sqrt{\alpha^2 + x^2} - 2\alpha  < x^2/\alpha
\quad \mbox{and} \quad 2 \sqrt{\alpha^2 + x^2} - 2\alpha  \to x^2/\alpha \quad \mbox{as} \quad \alpha \to \infty.$$
Hence, by the maximum principle, if $x(0) = 0$ then $x(t)$ is bounded above by solutions to the SNB differential equation
\begin{eqnarray}
 & &   \frac{dy}{dt} = \frac{y^2}{\alpha} - \mu, 
\qquad y(0) = 0, \quad \mu = 1 - \epsilon t. \label{SNBS2}
\end{eqnarray}
As $x(t) < y(t)$ it follows that $x$ tips later than $y$, and by continuity, the tipping time for $x$ converges to that for $y$ as $\alpha \to \infty$.

\vspace{0.1in}

\noindent We next consider the value of $\mu \equiv \mu_{TP}$ at which the solution of \eqref{SNBS2} satisfies $y = K \gg 1$ (as described in Section \ref{tipdef}). It is well known \cite{cody} that this equation has the asymptotic solution
$$y(t)= \alpha^{2/3} \epsilon^{1/3} \mbox{Ai}'(\alpha^{-1/3} \epsilon^{-2/3} \mu)/\mbox{Ai} (\alpha^{-1/3} \epsilon^{-2/3}\mu), $$
where Ai$(z)$ is the usual Airy function. 
If $\alpha^{2/3} \epsilon^{1/3} \ll K $ so that 
$$\alpha \ll K^{3/2} \epsilon^{-1/2} \equiv \alpha_1,$$
then for $y(t) = K \gg 1$ we must have $\alpha^{-1/3}\epsilon^{-2/3}\mu$ close to the first zero $-c_0$ of the Airy function. This gives the 'usual' estimate \eqref{1.7} for the tipping value.
In contrast if $\alpha > \alpha_1$ then the tipping value does not occur close to the rescaled zero of the Airy. In this case, as $\alpha \to \infty$ the estimate (\ref{c41}) applies.

\vspace{0.1in}

\begin{figure}[!htbp]
	\centering
 \includegraphics[width=0.46\textwidth]{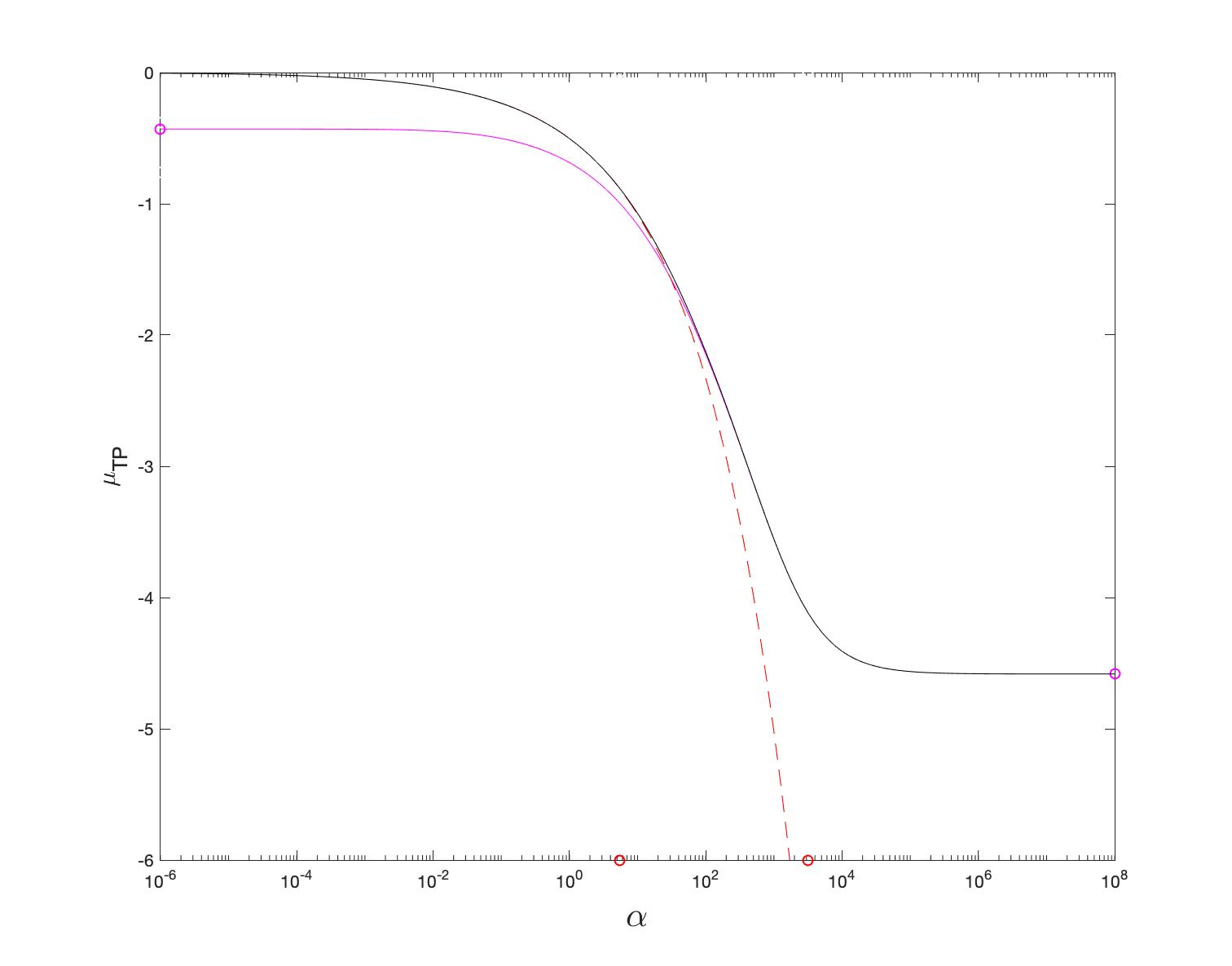}
		  \includegraphics[width=0.48\textwidth]{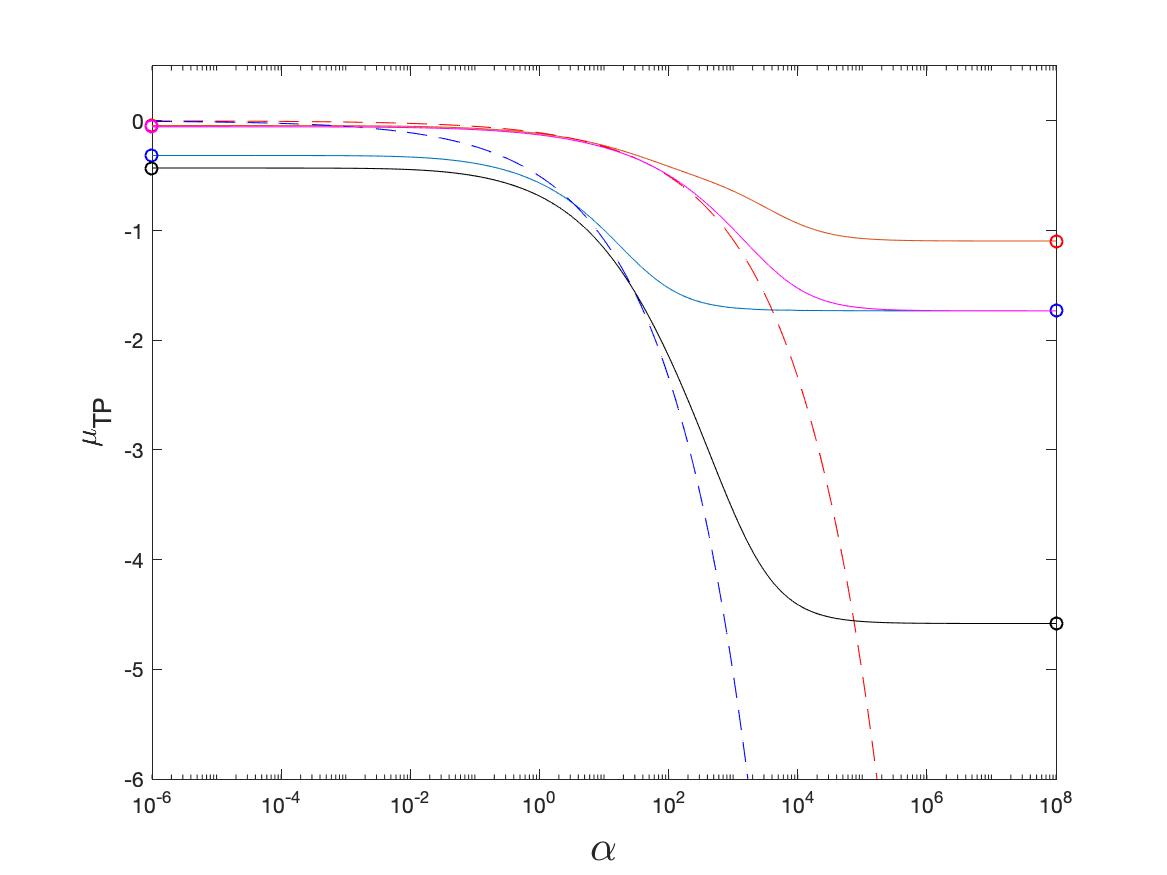}
	\caption{The tipping point $\mu_{TP}$ plotted as a function of $\alpha$ when $\mu = 1 - \epsilon t, x(0) = 0$. (left) we take $K = 100, \epsilon = 0.1$. In maroon is $\mu_{TP}$ and in black is the tipping point for (\ref{SNBS2}).  The $\alpha = 0$ and $\alpha = \infty$ estimates are given by circles, and the intermediate SNB estimate of $-c_0 \; \alpha^{1/3} \; \epsilon^{2/3}$ as a dashed line. The values of $\alpha_0, \alpha_1$ are given as circles on the x-axis. (right)  $K = 10$ and  $\epsilon = 0.01$ (red), $\epsilon = 0.1$ (blue), together with $K = 100$ and  $\epsilon = 0.01$ (maroon), $\epsilon = 0.1$ (black). Dashed lines again give the SNB estimate for the two values of $\epsilon$.}
	\label{fig:c2}
\end{figure}

\noindent The above calculations allow us to assess the impact of the drift rate, the smoothing and the limiting value $K$ on the unforced problem with slow drift.  We illustrate this in Figure \ref{fig:c2}, by considering the value of $\mu_{TP}$ as a function of $\alpha$  for $\epsilon = 0.1, 0.01$ and for  $K=10,100$.   
This figure also shows the asymptotic SNB prediction (\ref{1.7}) as dashed lines and the tipping point for (\ref{SNBS2}). We can clearly see the three asymptotic ranges for $\alpha$ with the SNB estimate reasonable for intermediate values of $\alpha_0 < \alpha < \alpha_1$.

\subsection{Smoothed oscillatory forcing}

\begin{figure}[htb!]
    \centering
    \includegraphics[width=0.5\textwidth]{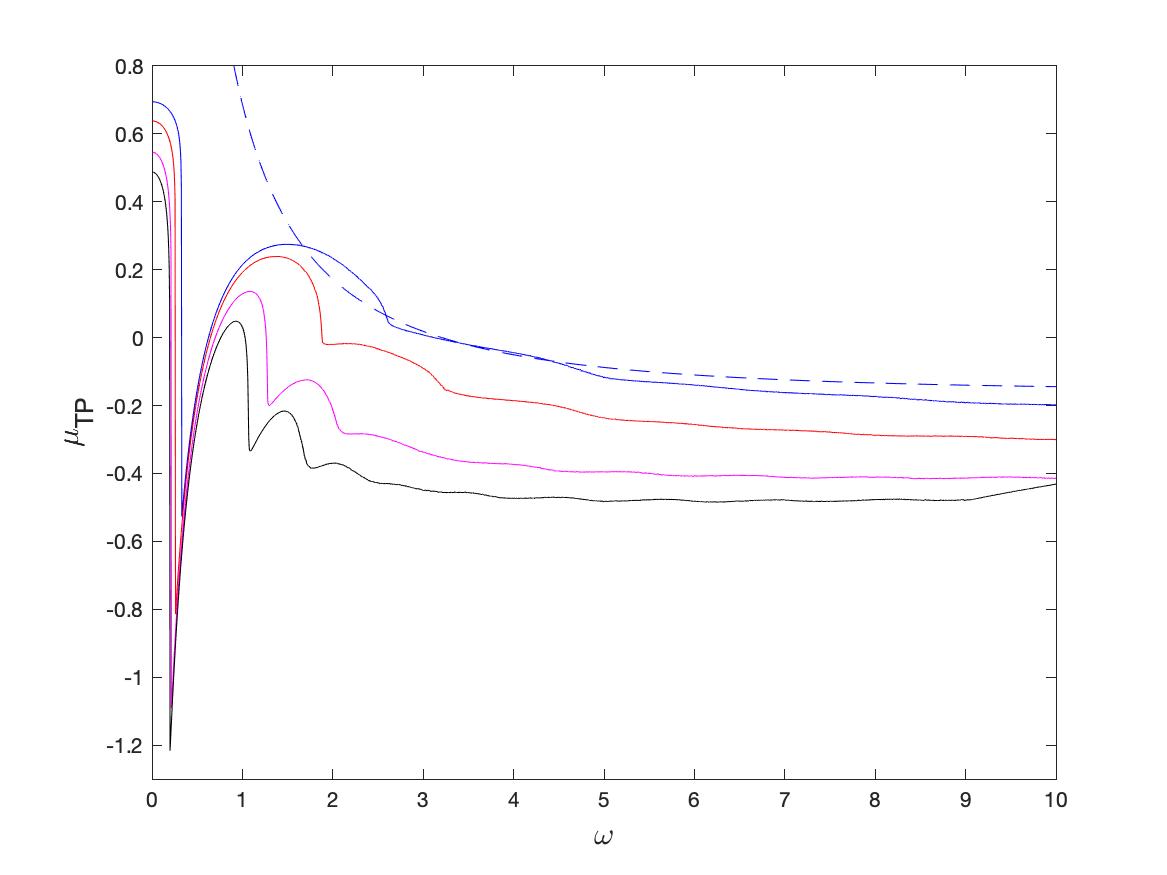}
    \caption{Tipping points $\mu_{TP}$ when $A = 1, \epsilon = 0.1, \omega \in [0,10]$ and smoothness $\alpha = 0,0.1,0.5,1$. There is a sharp transition at $\omega = 0.32$ in all cases, and evidence of further transitions for larger values of $\omega$. Dashed is the large $\omega$ approximation for $\alpha=0$ given in \eqref{3}. Note again that the value of $\mu_{TP}$ decreases as $\alpha$ increases. We compare this figure with those given in Section 4.}
    \label{fig:22ff}
\end{figure}

\noindent The algebraic method used to establish the location of the cyclic fold and related results, although very revealing, relies on the piece-wise linearity of the unsmoothed problem,
and cannot be applied directly to the smoothed system. However we may apply Lemma 6 to conclude directly that for any $A$ and $\epsilon > 0$ the value of $\mu_{TP}$ decreases with $\alpha$. Given that $\mu_{TP} \to \mu_{CF}$ as $\epsilon \to 0$ we may also deduce that $\mu_{CF}$ decreases with $\alpha$.
In Figure \ref{fig:22ff} we take $A = 1$ and plot $\mu_{TP} $ as a function of $\omega$ when $\epsilon = 0.1$ for different values of $\alpha$. We can see in this figure that the graph of $\mu_{TP}$ is  non-monotonic for smaller values of $\alpha$ with a sharp transition clearly evident in all cases. It is interesting that  all of the curves appear to have the transition point at $\omega \approx 0.32$ regardless of the value of $\alpha.$ This figure shows that the results for the location of $\mu_{TP}$ shown in the last section are robust to smoothing the original problem.

\section{Conclusions and future work}

\noindent We have studied the dynamics of tipping close to a non-smooth fold (NSF) in an oscillatory forced system  with slow drift. In this context we consider the influence of both the slow variation of a critical parameter $\mu$ and an external oscillatory forcing with amplitude $A$ and frequency $\omega$. Traditional studies of the detection of tipping in (for example) climate systems have centred around dynamic bifurcation near saddle-node bifurcations (SNB) in smooth systems. \\ 

\noindent In the SNB setting, with slow parameter drift, it is possible to make precise estimates of the location of the tipping points, and these estimates vary smoothly with the parameters in the system. These estimates are  typically made either by determining the system parameters, or by making measurements and observing a 'slowing down' in the system response as tipping is approached \cite{Lenton11}, thus identifying the lag of tipping relative to the related static SNB. Analyses for such  systems with an external oscillatory forcing show that an advance in the tipping combines additively with the lag  from the dynamic bifurcation, illustrated for both canonical and application-based models \cite{zhu2015tipping}.

\vspace{0.1in}

\noindent In this paper we exploit the linear structure for the one-dimensional reduction of the Stommel 2-box model, also appearing in other reduced climate models.
This provides the opportunity for explicit expressions for several important aspects of the system, not available in even the simplest nonlinear canonical models, such as the SNB.
Observations from our results, in comparison with SNB problem, indicate that predictions of tipping rely on a different balance of factors as compared with the smooth case, leading to the following conclusions:

\begin{enumerate}
\item Smoothing postpones tipping. In particular if predictions of tipping are made by modelling a problem with a saddle-node bifurcation, then the (possibly more representative) non-smooth, or lightly smoothed, problem can tip earlier. Or in other words, dramatic changes in behaviour effects may occur sooner than one might otherwise expect.

\item 
Tipping in non-smooth systems is not determined by the eigenvalues of the linearisation. In particular,
there is no equivalent of tipping occurring when an eigenvalue drops to zero. This rules out the identification of the closeness to tipping by monitoring the 'slowing down' in the behaviour associated with a zero eigenvalue.


Furthermore, the analysis of the reduced non-smooth model identifies  scenarios of significant uncertainty, relevant for both non-smooth and smooth systems: 

\item The critical value of the tipping parameter $\mu_{TP}$ does not behave monotonically for smaller values of $\omega$. Indeed we may see large transitions in the tipping times as parameters vary. This makes the estimation of tipping times in the context of noisy parameter values uncertain. We note that similar non-monotone behaviour for smaller values of $\omega$ can also be seen in the smoothed system and also in the forced SNB problem \cite{zhu2015tipping}.

\end{enumerate}

\vspace{0.1in}


\vspace{0.1in}

\noindent As a broader conclusion,  this analysis points to the need to  drawing conclusions  about the time and location of tipping points carefully when  based on  a smooth (unforced) saddle-node model. Nature generally has more complexities, such as disparate timescales and multiple contributing factors, which may motivate non-smooth problems as more realistic descriptions than the idealised smooth models studied in the literature.

\vspace{0.1in}

\noindent We shall explore these and additional features of non-smooth models  in forthcoming work that extends this analysis to the two-dimensional Stommel-Box model \cite{cody}, in which resonance effects play an important role. We will also look at the impact of additional stochastic forcing terms,  and the determination of tipping from noisy data gathered from a non-smooth system.

\bibliographystyle{plain}
\bibliography{algbib}

\begin{thebibliography}{10}

\bibitem{bernardo2008piecewise}
Mario Bernardo, Chris Budd, Alan~Richard Champneys, and Piotr Kowalczyk.
\newblock {\em Piecewise-smooth dynamical systems: theory and applications},
  volume 163.
\newblock Springer Science \& Business Media, 2008.

\bibitem{cody}
C.~Budd, C.~Griffith, and R.~Kuske.
\newblock Dynamic tipping in the non-smooth stommel-box model, with fast
  oscillatory forcing.
\newblock {\em Physica D:nonlinear phenomena}, 432, 2022.

\bibitem{dijkstra2013nonlinear}
Henk~A Dijkstra.
\newblock {\em Nonlinear climate dynamics}.
\newblock Cambridge University Press, 2013.

\bibitem{haberman79}
Richard Haberman.
\newblock Slowly varying jump and transition phenomena associated with
  algebraic bifurcation problems.
\newblock {\em SIAM Journal on Applied Mathematics}, 37(1):69--106, 1979.

\bibitem{kaper}
H.~Kaper and H.~Engler.
\newblock {\em Mathematics and Climate}.
\newblock SIAM, 2013.

\bibitem{KimWang18}
Wang~K.W. Kim, J.
\newblock Predicting non-stationary and stochastic activation of saddle-node
  bifurcation in non-smooth dynamical systems.
\newblock {\em Nonlinear Dyn}, 93:251–258, 2018.

\bibitem{Lenton11}
Timothy Lenton.
\newblock Early warning of climate tipping points.
\newblock {\em Nature Climate Change}, 1:201–209, 2011.

\bibitem{MoruBudd20}
K.S. Morupisi and C.~Budd.
\newblock {An analysis of the periodically forced PP04 climate model, using the
  theory of non-smooth dynamical systems}.
\newblock {\em IMA Journal of Applied Mathematics}, 86:76--120, 2021.

\bibitem{paillard2004antarctic}
Didier Paillard and Fr{\'e}d{\'e}ric Parrenin.
\newblock The antarctic ice sheet and the triggering of deglaciations.
\newblock {\em Earth and Planetary Science Letters}, 227(3-4):263--271, 2004.

\bibitem{stommel1961thermohaline}
Henry Stommel.
\newblock Thermohaline convection with two stable regimes of flow.
\newblock {\em Tellus}, 13(2):224--230, 1961.

\bibitem{bury}
I.~Pavithran M. Scheffer T. Lenton M.~Anand T.~Bury, R.~Smith and C.~Bauch.
\newblock Deep learning for early warning signals of tipping points.
\newblock {\em PNAS}, 118:1--9, 2021.

\bibitem{zhu2015tipping}
Jielin Zhu, Rachel Kuske, and Thomas Erneux.
\newblock Tipping points near a delayed saddle node bifurcation with periodic
  forcing.
\newblock {\em SIAM Journal on Applied Dynamical Systems}, 14(4):2030--2068,
  2015.

\end{thebibliography}

\section{Appendices}

\appendix

\section{Proof of Lemma 1}

\begin{proof} Trivially we have
\begin{equation}
    \dot{x} = -\mu -2x, \quad \mbox{if} \quad x < 0,
    \quad \mbox{and} \quad
   \dot{x} = -\mu +2x, \quad \mbox{if} \quad x > 0.
    \label{ca2}
\end{equation}

\noindent Hence, if $\mu$ is fixed, then the fixed points at 
$x^* = \mp \mu/2$ are {\em stable} if $x < 0$ and {\em unstable} if $x > 0.$

\vspace{0.1in}

\noindent Clearly we have $\mu = \mu_0 - \epsilon t.$
Both equations in (\ref{ca2}) then have exact solutions. The choice of initial condition implies that initially the system satisfies $x < 0$. A direct calculation then gives:
\begin{equation}
    x(t) = \left( \frac{\mu_0}{2} + \frac{\epsilon}{4} \right) e^{-2t} - \frac{\mu_0}{2} + 
    \frac{\epsilon t}{2} - \frac{\epsilon}{4} \equiv \left( \frac{\mu_0}{2} + \frac{\epsilon}{4} \right)  e^{-2t} - \mu/2 - \frac{\epsilon}{4}. 
    \label{ca6}
\end{equation}

\noindent If $\epsilon$ is {\em small} then $\mu = 0$ when $t = 1/\epsilon \gg 1$. Hence if $\mu$ is close to zero then terms of the form $e^{-2t} = {\cal O}(e^{-1/\epsilon})$ and can be neglected. Accordingly ignoring exponentially small terms we have:
$$x = -\mu/2 - \epsilon/4.$$
Hence we intersect the set $x = 0$ when 
$ \mu = -\epsilon/2.$ 
We have
$\dot{x} = -\mu = \epsilon/2 > 0$
so that the trajectory crosses from $x <0 $ to $x > 0$.
Assume that $x=0$ at time $t_0$. Set $s = t - t_0.$ Expressing derivatives with respect to $s$, (\ref{ca2}) then becomes:
$$\dot{x} = 2x - \mu, \quad \mbox{with} \quad
x(0) = 0, \quad \mu(s) = -\epsilon/2 - \epsilon s.$$
Again we can solve this system directly to give
\begin{equation}
    x(s) = \frac{\epsilon}{2} e^{2s} + \mu/2 - \epsilon/4.
\end{equation}
\noindent The system has tipped if $x(s) = K \gg 1$. We then have
$$\epsilon e^{2s}/2 - \epsilon s/2 - \epsilon/2 = K,
\quad \mbox{so that} \quad e^{2s} = \frac{2K}{\epsilon} + s + 1.$$
Thus $s$ satisfies the fixed point equation
$$s = \frac{1}{2} \log \left( \frac{2K}{\epsilon} \right) 
+ \frac{1}{2} \log \left(1 + \frac{\epsilon s}{2K} + \frac{\epsilon}{2K} \right).$$
Iterating this gives 
$$s  = \log (2K/\epsilon)/2 + \frac{\epsilon}{8K} \log \left( \frac{2K}{\epsilon} \right) + \ldots $$
which then yields \eqref{ca7}.
\end{proof}

\section{Proof of Lemmas 2(i),3, and 4}

\noindent By looking at the $\omega \gg 1$ limit we prove Lemma 2(i), Lemma 3 and Lemma 4.

\vspace{0.1in}

\noindent We start by proving Lemma 4. Motivated by (\ref{n10}) we note that if $A = {\cal O}(1)$ for large $\omega$ we have $\mu = {\cal O}(1/\omega).$ Accordingly we make a rescaling of the algebraic system above, setting
$$\mu = \nu/\omega, \quad C^{\pm} = D^{\pm}/\omega$$
and multiplying throughout by $\omega.$ We next expand all of the rescaled expressions in powers of $1/\omega$ up 
to ${\cal O}(1/\omega^2)$ and look at successive terms. 
This gives:
\begin{equation}
\frac{\nu}{2} + D^+ + A \sin(a)  - \frac{2 A \cos(a)}{\omega} - \frac{4 A \sin(a)}{\omega^2} =  {\cal O}(1/\omega^3) 
\label{31}
\end{equation}
\begin{equation}
\frac{\nu}{2} + D^+\left(1 + \frac{2(b-a)}{\omega} + \frac{2(b-a)^2}{\omega^2} \right) + A \sin(b)  - \frac{2 A \cos(b)}{\omega} - \frac{4 A \sin(b)}{\omega^2} = {\cal O}(1/\omega^3),
\label{32}
\end{equation}
\begin{equation}
 -\frac{\nu}{2} + D^- + A \sin(b)  + \frac{2 A \cos(b)}{\omega} - \frac{4 A \sin(b)}{\omega^2} = {\cal O}(1/\omega^3),
 \label{33}
 \end{equation}
 \begin{equation}
-\frac{\nu}{2} + D^-\left(1 - \frac{2(2\pi + a -b)}{\omega} + \frac{2(2\pi+a-b)^2}{\omega^2} \right) + A \sin(a)  + \frac{2 A \cos(a)}{\omega } - \frac{4 A \sin(a)}{\omega^2} = {\cal O}(1/\omega^3).
\label{34}
\end{equation}


\vspace{0.1in}

\noindent A first estimate gives:

$$\frac{\nu}{2} + D^+ +  A \sin(a)   = {\cal O}(1/\omega), \quad 
\frac{\nu}{2} + D^+ +  A \sin(b)  = {\cal O}(1/\omega),$$
$$-\frac{\nu}{2} + D^- + A \sin(b) = {\cal O}(1/\omega), \quad 
0 = -\frac{\nu}{2} + D^- + A \sin(a) = {\cal O}(1/\omega).$$

\noindent Thus $\sin(a) = \sin(b) + {\cal O}(1/\omega)$ so that there is a constant $\Delta$ with $b = \pi - a + \Delta/\omega + {\cal O}(1/\omega^2),$ so that 
\begin{equation}
 b - a = \pi - 2a + \Delta/\omega + {\cal O}(1/\omega^2), \quad 2\pi + a - b = \pi + 2a - \Delta/\omega  + {\cal O}(1/\omega^2),
\label{21}
\end{equation}
and hence
\begin{equation}
    \cos(a) +  \cos(b) = {\cal O}(\sin(a)/\omega).
    \label{21a}
\end{equation}
Also we deduce that there is a constant $D$ so that 
\begin{equation}
    D^- - D^+ = \nu + D/\omega + {\cal O}(1/\omega^2), \quad D^- + D^+ = -2 A \sin(a) + {\cal O}(1/\omega).
    \label{22}
\end{equation}

\vspace{0.1in}

\noindent Consider now the expressions given by (\ref{31})  - (\ref{34}) combined with the results in (\ref{21}).  After some manipulation this gives the identity
\begin{equation}
\nu + (D^+-D^-) + D^- \left( \frac{2(\pi + 2a)} {\omega} - \frac{2 \Delta}{\omega^2} - \frac{2(\pi +2a - \Delta/\omega)^2 }{\omega^2} \right) -\frac{4A \cos(a)}{\omega} = {\cal O}(1/\omega^3).
\label{k1}
\end{equation}
Similarly consider (\ref{33})-(\ref{32}). This gives
\begin{equation}
-\nu + (D^- - D^+) -  D^+ \left(  \frac{2(\pi - 2a)}{\omega} + \frac{2\Delta}{\omega^2} + \frac{2 (\pi -2a + \Delta/\omega)^2}{\omega^2} \right) + \frac{4A \cos(b)}{\omega} = {\cal O}(1/\omega^3).
\label{k2}
\end{equation}
Adding (\ref{k1}) and (\ref{k2}), multiplying by $\omega$ and applying (\ref{21}) again gives:
$$2\pi (D^- - D^+)  + 4a (D^- + D^+) - (D^+ + D^-) \frac{\Delta}{\omega}
-2(D^- + D^+) \frac{\pi^2}{\omega} - 8A \; \cos(a) = {\cal O}\left(\frac{a}{\omega},\frac{1}{\omega^2} \right).$$

\noindent Now, substituting (\ref{22}) we get 
$$2\pi \nu  - 8A \; (a \sin(a) + \cos(a)) = -2\pi \frac{D}{\omega} + {\cal O}\left(\frac{a}{\omega},\frac{1}{\omega^2} \right).$$

\vspace{0.1in}

\noindent If we now calculate (\ref{k1}) - (\ref{k2})  we have

$$2\nu + 2(D^+ - D^-) + \frac{1}{\omega} \left[ 2 \pi (D^+ + D^-) + 4 a (D^- - D^+) \right] + {\cal O}(1/\omega^2). $$
Hence
\begin{equation}
D^- - D^+ = \nu + \frac{2}{\omega} \left[ -\pi \sin(a) + a \nu \right] + {\cal O}(1/\omega^2).
\label{35}
\end{equation}
We see from this that in (\ref{22}) we have $D = 2(-\pi \sin(a) + \nu a) = {\cal O}(a).$ Thus we have
\begin{equation}
\pi \nu - 4 A \; (a \sin(a) + \cos(a)) = {\cal O}\left(\frac{a}{\omega},\frac{1}{\omega^2} \right).
\label{36}
\end{equation}
Setting now $\mu = \nu/\omega$ and dividing by $\omega$ gives (\ref{coct1}), which completes the proof of Lemma 4.

\qed 

\vspace{0.2in}

\noindent {\em Proof of Lemma 2(i)} To prove Lemma 2(i) we will study the form of the {\em cyclic fold bifurcation} arising from the formula in Lemma 4. To do this we will assume that $a$ is small. Substituting into (\ref{36}) we have
$$\frac{\pi \nu}{4} = A \left (1 +\frac{a^2}{2} \right) + {\cal O} \left(\frac{a}{\omega},\frac{1}{\omega^2} \right).$$
We immediately see that  there is a fold bifurcation when to leading order
$$a = 0, \quad \nu = \frac{4 A}{\pi} + {\cal O} \left( \frac{1}{\omega^2} \right).$$
To investigate this more precisely, we assume that the remainder terms in the expression (\ref{36}) are given by $P a/\omega + Q/\omega^2$
so that
$$\frac{\pi \mu}{4} = A \left(1 +\frac{a^2}{2} \right) + \frac{P a}{\omega} + \frac{Q}{\omega^2}.$$
If we assume that close to the cyclic fold we have $a = a_1/\omega$ so that $a = {\cal O}(1)$ then this expression balances and we have
$$\frac{\pi \mu}{4 \lambda} = 1 +\frac{1}{\omega^2} \left( a_1^2 + P a_1 + Q \right).$$
Then the cyclic-fold bifurcation occurs when 
$$a_1 = -P/2, \quad \mu = \frac{4 \lambda}{\pi} + \frac{-P^2/2 + Q}{\omega^2}.$$ 
The concludes the proof of Lemma 2 (i). 
\qed

\vspace{0.2in}

\noindent To prove Lemma 3 we collect up the results of the calculations above and rescale. This gives the estimates for $C^{\pm}$  and $a$. To estimate the average $\langle x \rangle$ of $x(t)$ we then integrate $x(t)$ directly over the interval $[0,2\pi/\omega]$ and divide by $2\pi/\omega$. Substituting the above estimates into the result gives (after some manipulation) the estimate for $\langle x \rangle$.  

\qed  

\section{Numerical estimates of the cyclic fold $\mu_{CF}$ when $A = 1$}

\vspace{0.2in}

\begin{center}
\begin{tabular}{|c|c|c|}
\hline  
$\omega$    & $\mu_{CF}$  & $\mu_G$  \\
\hline   
20 & 0.0635  & 0.0995    \\
15 & 0.0846  & 0.1322    \\
10 & 0.1265  & 0.19612   \\
8  & 0.15737  & 0.2425   \\
5 & 0.2471   & 0.3714    \\
4 & 0.3037  & 0.4472     \\
3 & 0.3933    & 0.5547   \\
2.5 & 0.4589 & 0.6247    \\
2 & 0.545 & 0.7071       \\
1.9 & 0.5661 & 0.7250    \\
1.8 & 0.5872 & 0.7433    \\
1.7 & 0.6101 & 0.7619    \\
1.6 & 0.635 &  0.7809    \\
1.5 & 0.666 & 0.8        \\
1.4 & 0.689 & 0.8192     \\
1.2 & 0.7458 & 0.8575    \\
1   & 0.804  & 0.8944    \\
0.8 & 0.875 & 0.9285     \\
0.5 & 0.945 & 0.9701     \\
0.3 & 0.9780  & 0.9889   \\
0.2 & 0.9902 & 0.995     \\
0.1 & 0.9975 & 0.9988    \\
0  & 1.000   & 1.000     \\
\hline
\end{tabular}
\end{center}

\end{document}